\documentclass[12pt]{amsart}
\usepackage{diagrams}
\makeatletter
\@mparswitchfalse
\def\page#1{}
\def\eqalign#1{\null\,\vcenter{\openup\jot\m@th
  \ialign{\strut\hfil$\displaystyle{##}$&$\displaystyle{{}##}$\hfil
      \crcr#1\crcr}}\,}
\def\eqalignbot#1{\null\,\vbox{\openup\jot\m@th
  \ialign{\strut\hfil$\displaystyle{##}$&$\displaystyle{{}##}$\hfil
      \crcr#1\crcr}}\,}
\let\atop\@@atop
\makeatother

\usepackage{amssymb}

\begin{document}

\title{Topology of Algebraic Varieties}

\author{Fouad Elzein}
\address{\hskip-\parindent
        Fouad Elzein\\
        D\'epart. de Math\\
        Univ. de Nantes\\
        44072 Nantes Cedex 03\\
        FRANCE}
\email{elzein@math.univ-nantes.fr}

\author{Andr\'as N\'emethi}
\address{\hskip-\parindent
        Andr\'as N\'emethi\\
        100 Math Building\\
        Ohio State University\\
        231 West 18th. Avenue\\
        Columbus, Ohio, USA 43210-1174}
\email{nemethi@math.ohio-state.edu}

\thanks{During the preparation of this work, the second author visited
Univ. de Nantes in June 1998 and the first author visited MSRI in
Berkeley in September 1998.  We thank these institutions for their
hospitality.
Research at MSRI is supported in part by NSF grant DMS-9701755.
The second author was partially supported by NSF grant DMS-9622724}

\begin{abstract}
Let $Y$ be a normal crossing divisor in the smooth projective algebraic variety
 $X$ (defined over ${\mathbb C}$)
and let   $U$ be  a tubular neighbourhood of $Y$ in $X$.
We construct homological cycles generating   
$H_*(A,B)$, where $(A,B)$ is one of the following pairs
$(Y,\emptyset)$,\ $(X,Y)$,\ $(X,X-Y)$,\ $(X-Y,\emptyset)$ and 
$(\partial U,\emptyset)$.
The construction is compatible with the 
weights in $H_*(A,B,{\mathbb Q})$ of  Deligne's mixed Hodge structure.
\end{abstract}

\def\thesection{\Roman{section}}

\def\fighere{\par$$\mbox{FIGURE GOES HERE}$$\par}
\let\d\partial
\def\EE{\mathcal E}
\def\C{\mathbb C}
\def\Q{\mathbb Q}
\def\R{\mathbb R}
\def\Z{\mathbb Z}
\def\Res{\operatorname{Res}}
\def\im{\operatorname{im}}
\def\Ker{\operatorname{Ker}}
\let\ker\Ker
\def\Tot{\operatorname{Tot}}
\def\Gr{\operatorname{Gr}}
\def\Hom{\operatorname{Hom}}
\def\Log{\operatorname{Hom}}
\def\codim{\operatorname{codim}}
\def\dim{\operatorname{dim}}
\let\emptyset\varnothing
\def\resp{\operatorname{respectively}}
\def\resp{\operatorname{respectively}}
\def\mod{\operatorname{mod}}
\def\modulo{\operatorname{modulo}}
\def\b{\bullet}
\def\c{\mathcal}
\def\p{\pitchfork}
\def\Y{Y_{\alpha_1,\ldots,\alpha_p}}
\def\u{\underline}
\def\y{\tilde{Y}^p}
\def\s{\star}

\newenvironment{theorem}
        {\noindent\bgroup{\bf Theorem.} \it\ignorespaces}
        {\egroup \medbreak}
\newenvironment{prop}
        {\noindent\bgroup{\bf Proposition.} \it\ignorespaces}
        {\egroup \medbreak}
\newenvironment{lemma}
        {\noindent\bgroup{\bf Lemma.} \it\ignorespaces}
        {\egroup \medbreak}
\newenvironment{corollary}
        {\noindent\bgroup{\bf Corollary.} \it\ignorespaces}
        {\egroup \medbreak}
\newenvironment{remark}
        {\noindent\bgroup{\bf Remark.} \ignorespaces}
        {\egroup \medbreak}

\maketitle

\section{Introduction}

In his fundamental article [L], Leray gives a cohomological version of
the theory of residues  along a smooth hypersurface in an analytic
manifold.  The original motivation of the present article was 
the generalization of the theory of residues, but now along a normal
crossing divisior ({\em NCD}) $Y$ of  a smooth projective  algebraic
variety $X$ (defined over $\C$). 
 Leray's construction of the inverse image of a cycle on $ Y$ 
leads to a homological cycle of the 
boundary  of a tubular neihgbourhood  of $Y$.
In our case, 
in  order to obtain a cycle as inverse image, we need to
impose
specific condition on the cycle on $ Y$. We obtain such conditions by
duality from 
Deligne's mixed Hodge structure  ({\em MHS})  on the cohomology of
algebraic
varieties [D]. 

Once this fact has been well understood, we realised that we needed
to develop a complete homological theory
of some complexes (compatible with the natural mixed Hodge structures),
and to describe  the topology/homology
of the pairs $(Y,\emptyset)$,\ $(X,Y)$,\ $(X,X-Y)$,\ $(X-Y,\emptyset)$ and 
$(\partial U,\emptyset)$, where $U$ is  a tubular neihgbourhood of $Y$ in $X$.
These results are presented in this paper. The generalization of the 
residue theorems will be published elsewhere. 

The reader may have recognized that the reduction to the {\em NCD} case uses
Hironaka's desingularization theorem [H]
in order to carry our results to smooth open algebraic varieties.  
In the cohomology theory,
the main ingredient is Deligne's mixed Hodge structure. 
By duality, one can define the weight filtration in homology with 
$\Q$-coefficients.
Our ultimate result is the construction of representative homological cycles
according to their weight.  Such constructions have been used
occasionally in examples or some specific situations.  Our objective
is a systematic treatment giving a general account of the topology of
algebraic varieties.

>From the point of view of the 
present article, the key argument  is the degeneration at rank two of
the weight spectral sequence defined by the mixed Hodge  complex,
valid only for algebraic varieties.
This in turn leads us to consider  homological complexes with good 
intersection theory and for this it has been necessary
to use the simplicial theory provided 
by the triangulations of varieties  instead of singular chains.

The paper is organized as follows: in the introduction we recall some
known facts about the topology of a regular neighbourhood of a {\em NCD}.  We
refer to the articles by Clemens [C] and A'Campo [A]. 
In the second part of this introduction,
 we already present the main steps of our construction.

In the second section we consider the  complex of ``dimensionally
transverse'' sub-analytic  chains. This is in fact
the complex of the sub-analytic  ``intersection chains''
  associated with 
the natural stratification of $Y$ and zero perversity
 (in the sense of Goresky--MacPherson [GM]). Moreover, we define 
and analyze different chain morphisms whose role is absolutely 
fundamental in the construction of cycles. 

Using the usual and the ``dimensionally transverse'' chains, in section 3
we 
define several homological double complexes $A_{**}(A,B)$
corresponding to the above pairs of spaces $(A,B)$. 
They carry  natural weight filtrations,
giving rise to the corresponding  spectral sequences. For example,
in the case of $(X-Y,\emptyset)$ 
it is possible to show that our spectral sequence
is dual to the weight spectral sequence of the logarithmic complex
via residues. In this section we prove for all of the considered pairs
 the degeneration property at rank 2 of the weight spectral sequence. 
In section 4 we review some of the  algebraic properties 
of a  double complex whose spectral sequence degenerates at rank 2.

Using the result of the previous two sections, in section 5 we construct
cycles compatible with the weight filtration. 
Here an additional important ingredient is the construction of 
chain morphism $Tot_*(A_{**}(A,B))\to C_*(A,B)$ for any pair 
$(A,B)$. This is based on the
results of section 3. (Here $Tot_*$ stays for the total complex, and $C_*(A,B)$
is the complex of sub-analytic chains of $(A,B)$.)

The case of $(\partial U,\emptyset)$ 
is the most involved, so we separated it in section 6.

\vspace{2mm}

\subsection{} Now, we start a more detailed  presentation.

\vspace{2mm}

\subsubsection{}{\bf The stratification of $Y$.} \  
Let $X$ be a smooth projective algebraic variety containing 
a normal crossing divisor $Y$.  Let $(Y_i)_{i\in{I}}$ be 
the decomposition of $Y$ into irreducible components. 
Here the index set $I$ is a totally ordered set. 

The divisor $Y$ has a natural stratification by subspaces $Y^j$
consisting of ``points of multiplicity $\geq j$ in $Y$'', i.e.:
 \[Y^j=
\bigcup_{\alpha_1<\cdots<\alpha_j} Y_{\alpha_1}\cap\cdots \cap Y_{\alpha_j},\
\ \mbox{with normalization} \ \ 
\tilde{Y}^j= \coprod_{\alpha_1<\cdots<\alpha_j} Y_{\alpha_1}
\cap\cdots\cap Y_{\alpha_j}.\] 
Sometimes we will use the notation $Y_{\alpha_1,\ldots,\alpha_j}=Y_{\alpha_1}
\cap \cdots\cap Y_{\alpha_j}$ as well.
It is convenient to write $\tilde{Y}^0=Y^0=X$.

\subsubsection{}{\bf The tubular neighbourhood of $Y$.}
For any variety  $X$  and divisor $Y$, $X$ admits a
triangulation  compatible with $Y$, hence there exists a regular (tubular)
neighbourhood in $X$ which is a deformation retract of $Y$ (see e.g. [R-S]). 
In the case of a {\em NCD},  the tubular neighbourhood $\, U\,$of $Y$ 
and a projection from $U$  to $Y$ can have even some additional 
properties, see [C].

A' Campo [A] describes such a  neighbourhood as follows. 
The real (oriented) blow-up
$\Pi{_i}: {Z_i \rightarrow X}$ with center  $Y_i$ provides a
differentiable manifold $Z_i$ with boundary and a projection $\Pi_i$
on $X$, inducing a diffeomorphism outside $Y_i$. Moreover,
$\Pi_i^{-1}(Y_i)$ is an $S^1$-bundle associated to the oriented
normal bundle $N_{Y_i/X}$.
We can identify the fibers at a point $y$ in $Y_i$ with the set of
real oriented normal directions to $Y_i$.  For example, if $Y_i$
is a point in $\C$, then $Z_i$ is a half-cylinder with boundary $S^1$
over the point.  
In general, the boundary of $Z_i$, equal to
$\Pi_i^{-1} (Y_i)$, is diffeomorphic to the boundary of the complement
of an open tubular neighbourhood of $Y_i$ in $X$.

Next we consider the fibered product of the projections $(\Pi_i)_{i\in
I}$ over $X$.  Thus we obtain $\Pi : Z\rightarrow X$, where $Z$ is a
manifold with corners 
 whose  boundary $\partial Z$ equal to $\Pi^{-1} (Y)$. 

Here $Z$ is
homeomorphic to the complement of an open tubular neighbourhood $U$ of
$Y$ in $X$, hence $\partial Z$ is homeomorphic also to the boundary 
of the closure of $ U$.
 In fact $U$ is homeomorphic to the mapping
cylinder of $ {\Pi}|{\partial Z}: \partial Z \rightarrow Y$.

In order to see this, let's fix a homeomorphism 
$\phi_{\alpha}:\partial Z\times [0,\alpha]\to Z_{\alpha}$, where 
$Z_{\alpha}$ is some collar, and 
$\phi_{\alpha}|\partial Z\times \{0\}$ is the inclusion 
$\partial Z\hookrightarrow Z_{\alpha}$.
There is a natural projection $p_{\alpha}:Z_{\alpha}\to \partial Z$ given by
$pr_1\circ \phi^{-1}_{\alpha}$, where $pr_1$ is the projection on the 
first factor. 
(This retract can be extended to a strong deformation retract of the 
inclusion $\partial Z\hookrightarrow Z_{\alpha}$.)
Using this, one can define  tubular neighbourhoods
of  $Y$ in $X$ as follows.

For any $\epsilon \in (0,\alpha]$, define 
$U_{\epsilon}:=\Pi(\phi_{\alpha}(\partial Z\times [0,\epsilon]))$. 
Its boundary $\partial U_{\epsilon}$ is clearly  
$\Pi(\phi_{\alpha}(\partial Z\times \{\epsilon\}))$. Moreover, there is a unique map 
$q_{\epsilon}:U_{\epsilon}\to Y$ such that $q_{\epsilon}\circ  \Pi=\Pi\circ
p_{\epsilon}$. From the definitions follows that for two neighbourhoods
$U_{\epsilon_1}\stackrel{i}{\hookrightarrow}U_{\epsilon_2}$, for 
 $0<\epsilon_1< \epsilon_2$, one has $q_{\epsilon_2}\circ i=q_{\epsilon_1}$.
This shows, that one can construct a uniform family of tubular 
neighbourhoods $U_{\epsilon}$ of $Y$ in $X$,  all of them homeomorphic
to the mapping cylinder of $\Pi|\partial Z$, 
with all the natural compatibilities.
This construction is rather satisfactory in the discussion of topological
invariants.

But, in fact, $Z$ has a natural semi-analytic  structure (see II.3.1),
and the above homeomorphism $\phi_{\alpha}$ can be chosen as a sub-analytic 
homeomorphism as well (in fact, even
as a piecewise analytic isomorphism with respect to the natural
stratification). For the proof of this last statement,
 see [P]. Therefore, all the homeomorphisms, discussed in the previous
paragraphs, can be considered as sub-analytic homeomorphisms. 

\subsubsection{} {\bf The construction of the double complexes $A_{**}(A,B)$.}
For any $p\geq 0$, we denote the
complex of the usual (respectively the ``dimensionally transverse'')
sub-analytic  chains (cf. II.1) defined on $\y$ by $(C_*(\y),\partial)$
(respectively by 
$(C_*^\p(\y),\partial)$).
One has the following chain morphisms:

$$\begin{array}{ccccccccccc}
\cdots&
C_{*-4}^{\p}(\tilde{Y}^2)&\stackrel{\cap}{\leftarrow}&
C_{*-2}^{\p}(\tilde{Y}^1)&\stackrel{\cap}{\leftarrow}&
C_{*}^{\p}(\tilde{Y}^0)& &&&& \\
&&&&& \Big\downarrow\vcenter{\rlap{$j_1$}}&&&&& \\
&&&&&
C_*(\tilde{Y}^0)& \stackrel{i}{\leftarrow}&
C_*(\tilde{Y}^1)& \stackrel{i}{\leftarrow}&
C_*(\tilde{Y}^2)& \cdots \end{array}$$

\noindent Then define (modulo some shift of indexes):
$A_{**}(X,X-Y)=(C_*^\p(\y),\partial, \cap)_{p\geq 1}$,
$A_{**}(X-Y)=(C_*^\p(\y),\partial, \cap)_{p\geq 0}$,
$A_{**}(Y)=(C_*(\y),\partial,i)_{p\geq 1}$, 
$A_{**}(X,Y)=(C_*(\y),\partial,i)_{p\geq 0}$, and finally
$A_{**}(\partial U)$ by the cone of the  morphism $j_1$ in the above 
diagram where we replace $\tilde{Y}^0$ by $U$.  
The weight filtration $W_*$ of a double complex $A_{**}$ 
is defined  by $W_s:=\oplus_{p\leq s}A_{pq}$.

\subsubsection{}{\bf The construction of the cycles.}
First, using the geometry of the pair $(X,Y)$ (e.g. the projection $\Pi$
and intersections  corresponding to the stratification),
for any pair $(A,B)$ we construct a quasi-isomorphism
$m_{A,B}:Tot_*(A_{**}(A,B))\to C_*(A,B)$. Then
any representative $c_{st}$ of an element $[c_{st}]\in ker (d^1| E^1_{st})$ 
 is completed to a chain $c^{\infty}_{st}=c_{st}+c_{s-1,t+1}+\cdots
\in Z^{\infty}_{st}$ with $D(c^{\infty}_{st})=0$, where $D$ is the 
differential of the associated total complex $Tot_*(A_{**})$
(and $c_{pq}\in A_{pq}$). Therefore, $m_{A,B}$ associates with 
any $[c_{st}]$ a closed cycle $m_{A,B}(c^{\infty}_{s,t})$ of dimension $k=s+t$.

For the various pairs $(A,B)$ considered above, we have the following result.

\subsubsection{}{\bf Theorem.} {\em

a)\ The first term $(E^1_{st},d^1)$ of the weight spectral sequence
can be explicitly determined from
the homology of the spaces $\y$ and from the various normal bundles of the 
components of $\y$ in $\tilde{Y}^{p-1}$.

b)\ $E^r_{st}\Longrightarrow H_{s+t}(A,B,\Z)$, and  induces a weight filtration
on the integer homology. Moreover $E^{\infty}_{st}\otimes \Q=
Gr^W_{-t}H_{s+t}(A,B,\Q)$ (the last considered in Deligne's MHS).

c)\ $E_{**}^*\otimes \Q$ degenerates at rank 2.

d)\ The above construction provides all the cycles of the pair $(A,B)$
(modulo boundary) according to their weights. More precisely,
the homology class $[m_{A,B}(c^{\infty}_{s,t})]$ of dimension $k=s+t$
is well-defined modulo $W_{-t-1}H_k(A,B,\Q)$, and these type of classes
generate $W_{-t}H_k(A,B,\Q)$. }

\vspace{2mm}

\noindent 
Here we refer the reader to theorem III.0 and propositions V.2.2,\,  VI.2.3 \,
and \, VI.2.4\, for precise statements. 

\subsubsection{}{\bf Example.}\  Suppose
$Y=Y_1\cup Y_2$ has two components with smooth intersection $Y_{1,2}=Y_1\cap
Y_2.$ Two homology  classes $[a_i]\in H_k(Y_i),\  i=1,2$,  satisfying
$[a_1]\cap[Y_{1,2}]=[a_2]\cap[Y_{1,2}]$ in $H_{k-2}(Y_{1,2})$ give rise
in $\partial U$ to a $k+1$ dimensional cycle  (generating a homology class
of weight $-k-2$).  Indeed, first we can assume  that the 
representative $a_i$ is transversal to
$Y_{1,2}$ in $Y_i$.  Due to the condition about  $[a_i]\cap[Y_{1,2}]$, there
exists a chain $a_{1,2}$ in $Y_{1,2}$ such that $\partial
a_{1,2}=a_2\cap Y_{1,2}-a_1\cap Y_{1,2}$.  Since  $\dim \, a_i=k$, we have:
 \[\dim \, a_i\cap Y_{1,2}=k-2,
\
\dim a_{1,2}=k-1,\ \dim \Pi^{-1} a_i=k+1,\ \dim \Pi^{-1} a_{1,2}=k+1.\]
Moreover, by the very construction, $\Pi^{-1}(a_1+a_2+a_{1,2})$
has no boundary. It is the wanted closed cycle in $\partial U$.

\vspace{1mm}

Now, consider the case when $Y$ has three irreducible components
$Y_i$ ($i=1,2,3$) with $Y_{1,2,3}\not=\emptyset$. Similarly as above,
we want to lift some closed cycles  $a_i$ into $\partial U$. 
The first obstruction
is $[a_i\cap Y_{i,j}]=[a_j\cap Y_{i,j}]$ in the homology of $Y_{i,j}$ (for any
pair $i\not=j$). Using this, we create the new chains $a_{i,j}$ in $Y_{i,j}$ as
above. Now, if we want to lift these new chains and glue them together,
we face the second obstruction provided by the triple intersection
$Y_{1,2,3}$. The main point is that {\em 
there exists a good choice of $a_{i,j}$
such that this new obstruction is trivial}; its
triviality is equivalent with the vanishing of the second differential
$d_2$ of the spectral sequence. Basically, this is the main message of the
present paper.

\section{Topological Preliminaries}

\subsection{}{\bf Geometric chains.}

\subsubsection{} {\bf Preliminary remarks.} \ 
In some cases it is not absolutely evident how can we  dualize a result 
established in cohomology. For example, if we want to compute the 
homology of $Y$, then Deligne spectral sequence in cohomology has a very 
natural analog in homology. On the other hand, if we want   to find the
homological analogue of the cohomological theory of $X-Y$, then we have to
 realize that there is no obvious homological candidate.
 Actually, the $E_1$-term of the spectral sequence  associated with 
the  log complex $(\Omega_X^*(\log Y), W)$ can be easily dualized;
but we want (and need to) dualize the whole spectral  sequence, 
in particular we have to construct the $E^0$ term on $X$.  
This is crucial in the construction of 
cycles in $X-Y$ as well. 

Since  the {\em Gysin differential}, in the cohomological $E_1$-term of $X-Y$,
 dualizes to the 
intersection of cycles, we need to work with  {\em chains with good 
intersection properties}  with respect to the 
stratification defined by $Y$. 
Similarly as in the case of the intersection homology groups, we 
have several options to define our chain complex.  In the original definition 
of the  intersection homology groups, Goresky and MacPherson used 
geometric chains with some restrictions provided by the perversities [GM]. 
On the other hand, H.King recovered these groups using singular chains [K].
We will follow  here the first option (in fact, we will use sub-analytic
geometric chains). In order to have 
a good intersection theory, we need some kind of transversality 
property. Here again  we have several possibilities. Our choice
asks only a ``dimensional transversality'' of the chains. 
This has the big advantage that the complex of these chains
 coincides with the complex of zero perversity
chains of Goresky--MacPherson;  but has the disadvantage, that in the 
intersections of the cycles we have to handle an intersection multiplicity
problem. 

Finally, our choice for {\em sub-analytic} chains is motivated by the fact
that these chains are stable with respect to the real blowing up along $Y$
(in contrast with the P.L. chains). 

\subsubsection{} {\bf Good class of subsets and chains.}\
Before we start the precise definition of our  geometric chains, let us 
review briefly the general theory. We will follow the presentation from
[G], pages 146-155. 

Fix a manifold $M$. 
The group of {\em geometric chains} is defined in three 
steps. First, one defines {\em a good class ${\c C}$ of subsets of $M$}. 
This should satisfy the following properties:

(1)\ if a subset $S$ of $M$ is in ${\c C}$, then $M$ has a Whitney 
stratification such that $S$ is a union of strata, and each
 stratum is in ${\c C}$; 

(2)\ the class ${\c C}$ is closed under unions, intersections and
differences; 

(3)\ the closure of a subset in ${\c C}$ is in ${\c C}$. 

\vspace{1mm}

In the second step, one defines the {\em geometric prechains} (relative to 
${\c C}$). 
Basically, they can be represented as $\sum_{\alpha}m_{\alpha}S_{\alpha}$,
where each  $S_{\alpha}$ is a closed subset from ${\c C}$ (with a fixed
orientation) and $m_{\alpha}$ is an integer (the ``multiplicity of
$S_{\alpha}$''). A geometric chain is an equivalence class of geometric
prechains with respect to a natural equivalence relation. 
(For details, see [G].)

\vspace{1mm}

In this paper we will use sub-analytic chains associated with 
our  complex analytic manifolds: if $M$ is a complex analytic manifold, and 
we fix a real analytic structure on $M$, then the class of all
sub--analytic subsets of $M$ form the good class of subsets ${\c C}$.

\subsubsection{}{\bf The definition of the (geometric) chains.} \  
Now, we  define our chains. In some of the 
 definitions we follow closely the paper [GM] of Goresky and
MacPherson (where the case of the P.L. geometric chains is considered). 
We will regard our manifold $X$  as a pseudomanifold 
of dimension $2n$, where $n$ is the complex dimension of $X$.
We consider its natural stratification $X_{2n-2p+1}=X_{2n-2p}=Y^p$
for $p\geq 1$, hence $\Sigma=Y$. Similarly as above $\tilde{Y}^0=Y^0=X$.  

By [H2], $X$ admits a (canonical) sub--analytic triangulation, compatible with
the stratification. In fact, any two sub--analytic
triangulation admits a common refinement (see [H2] 2.4).

For  the next definitions,
 we fix an integer $p\geq 0$. Then $\tilde{Y}^p$ is again a pseudomanifold.
This stratification  will be denoted by 
$\{\tilde{Y}^p_{2n-2p-2r}\}_{r\geq 0}$.

If $T$ is a sub--analytic triangulation, 
let $C_*^T(\tilde{Y}^p)$ be the chain complex  of
simplicial chains of $\tilde{Y}^p$ with respect to $T$. By definition,
a chain of $\tilde{Y}^p$ is an element of $C_*^T(\tilde{Y}^p)$ for some 
sub--analytic triangulation
$T$, however one identifies two chains $c\in C_*^T$ and $c'\in C_*^{T'}$ if 
their canonical images in $C_*^{T''}$ coincide, for some common
refinement $T''$ of $T$ and $T'$.
The group of all chains is 
denoted  by $C_*(\tilde{Y}^p)$. Obviously, there is a natural boundary
operator $\partial$ which makes $C_*(\y)$  a complex. 

Notice that $C_*(\tilde{Y}^p)$ is exactly the group of sub-analytic 
geometric chains discussed in II.1.2. Indeed, if we fix an arbitrary 
closed 
sub-analytic subset $S$ of $\tilde{Y}^p$, then by [H2] there is a sub-analytic
triangulation $T$ which makes $S$ an element of $C_*^T(\tilde{Y}^p)$ 
(with all multiplicities one). 

If $\xi\in C_k^T$, 
then the support $|\xi|$ of $\xi$ is the union of the closures
of those $k$--simplices $\sigma$ for which the coefficient of $\xi$ is 
non--zero. Actually, the support of $\xi$ is independent on $T$, and
it is a $k$--dimensional sub-analytic  subset.  
We emphasize (even if this is clear in most of the 
 cases, since $X$ is compact), that in all our
discussions, we will deal with chains with {\em compact} supports. 
(See also  [B], I.) This remark is essential, when we 
will consider the sheaf-versions of our complexes.

The homology of the complex $C_*(\y)$ is the usual homology $H_*(\y)$ (cf.
[G], page 155). Moreover, for any triangulation $T$, the natural morphism 
$C^T_*(\y)\to C_*(\y)$ is a quasi-isomorphism. Indeed, the triangulation
homeomorphism $t:|K|\to \y$ from a simplicial complex $|K|$ to $\y$
identifies the simplicial homology of $|K|$ and the usual homology of $\y$. 

\subsection{}{\bf  ``Dimensionally transverse'' chains.}\\

We introduce the following double complex of 
{\em dimensionally transverse} chains in  $\tilde Y^p, p\geq 0$. 
We say that a chain $\xi\in C_k(\tilde{Y}^p)$ is {\em dimensionally 
transverse}, if 

\vspace{2mm}

a)\ $\dim(|\xi|)\cap \tilde{Y}^p_{2n-2p-2r}\leq k-2r$, for any $r>0$; and

b)\ $\dim(|\partial \xi|)\cap \tilde{Y}^p_{2n-2p-2r}\leq k-1-2r$, for any 
$r>0$. 

\vspace{2mm}

The subgroup of $C_*(\tilde{Y}^p)$ consisting of the {\em dimensionally 
transverse} chains  is denoted by  
$C_*^\pitchfork(\tilde Y^p)$. 
The boundary operator $\partial $ maps dimensionally transverse chains
in dimensionally transverse chains, hence defines a complex with 
$\partial^2=0$. 

The reader,  familiar with intersection homology, immediatelly will realize 
that the complex $C_*^{\pitchfork}(\tilde{Y}^p)$ is exactly the complex
$IC_*^{\bar{0}}\, (\tilde{Y}^p)$
of intersection chains corresponding to the zero perversity (see [GM]). 

\subsubsection { }
\begin{lemma} The natural inclusion $j:(C_*^{\p}(\tilde{Y}^p),\partial)
\to (C_*(\tilde{Y}^p),\partial)$ is a quasi-isomorphism. In particular, 
the homology of $(C_{*}^{\pitchfork}(\tilde Y^p),\partial)$ is
$H_{*}(\tilde Y^p)$.
\end{lemma}

\begin{proof} \ This follows from section 4.3 of [GM] (because
the Poincar\'e map $H^{2n-2p-k}(\tilde{Y}^p)$ $\to H_k(\tilde{Y}^p)$ is an
isomorphism, provided by the smoothness of $\tilde{Y}^p$). Actually,
since $\tilde{Y}^p$ is smooth, and the intersection homology group is 
independent on the  stratification (3.2 in [GM]), all the intersection
homology groups are the same and equal to the usual homology.
\end{proof}

\subsubsection{}{\bf The intersections with the strata.}\ We want 
to define the intersection $\xi\cap Y_{\alpha}\in 
C_{k-2}^{\pitchfork}(Y_{\alpha_1,\ldots,\alpha_p}$ $\cap Y_{\alpha})$  for any
chain $\xi=\xi_{\alpha_1<\cdots<\alpha_p}
\in C_k^{\pitchfork}(Y_{\alpha_1,\ldots,\alpha_p})$, where
$\alpha\not\in \{\alpha_1,\ldots,\alpha_p\}$.
(Here, if $p=0$ then $Y_{\alpha_1,\ldots,\alpha_p} $ denotes $X$.)
First notice the following fact (cf. [GM],  page 138).

\vspace{2mm}

\subsubsection{}{\bf Fact.} {\em 
\ If $C\subset \tilde{Y}^{p+1}$ is a ($k$--2)--dimensional sub-analytic
  subset of
$\tilde{Y}^{p+1}$ and if $D\subset C$ is a ($k$--3)--dimensional 
sub-analytic subset,
then there is a one-to-one correspondence  between chains $\beta\in
C_{k-2}(\tilde{Y}^{p+1})$ such that $|\beta|\subset C$,\ $|\partial\beta|
\subset  D$, 
and between homology classes $\tilde{\beta}\in H_{k-2}(C,D)$. Furthermore,
$\partial \beta$ corresponds to the class $\partial _*(\tilde{\beta})$ in
$H_{k-3}(D)$ under the connecting homomorphism $\partial_*:H_{k-2}(C,D)\to
H_{k-3}(D)$. }

\vspace{2mm}

In the definition of the intersection $\xi\cap Y_{\alpha}$, 
 we face the problem 
of the notion of an ``intersection multiplicity''. This is clarified 
in the work of Lefschetz, however we will work in the spirit of [GM]. 
It is clear that $|\xi|\cap Y_{\alpha}$ has dimension $\leq k-2$, and this 
intersection satisfies the transversality restrictions of the new
space. We have to determine the coefficients of the simplices which support
the intersection. 

Using the above Fact, the {\em intersection} $\xi\mapsto \xi\cap Y_{\alpha}$ is
completely determined by the following composition:

\vspace{2mm}

$$H_k(|\xi|,|\partial \xi|)$$
$$\hspace{2cm}\approx\ \uparrow \ \ \cap[Y_{\alpha_1,\ldots,\alpha_p}]$$
$$H^{2n-2p-k}(Y_{\alpha_1,\ldots,\alpha_p}-|\partial \xi|,
Y_{\alpha_1,\ldots,\alpha_p}-|\xi|)$$
$$\hspace{1cm} \downarrow \ \ i^*$$
$$H^{2n-2p-k}(Y_{\alpha_1,\ldots,\alpha_p}\cap Y_{\alpha}-|\partial 
\xi|,Y_{\alpha_1,\ldots,\alpha_p}\cap 
Y_{\alpha}-|\xi|)$$
$$\hspace{3cm}\approx\ \downarrow\ \ \cap [Y_{\alpha_1,\ldots,\alpha_p}
\cap Y_{\alpha}]$$
$$H_{k-2}(|\xi|\cap Y_{\alpha},|\partial \xi| \cap Y_{\alpha}).$$

\vspace{2mm}

\noindent Above, the first map is the (inverse) of the cap product
with the fundamental class of $Y_{\alpha_1,\ldots\alpha_p}$, 
as it is presented in the Appendix of [GM], page 162; the second map is
the restriction $i^*$, where $i$ is the natural inclusion;  and finally,
the third map is again 
the cap product with the fundamental class of $Y_{\alpha_1,\ldots
\alpha_p}\cap Y_{\alpha}$.
The cap products are isomorphisms (cf. [{\em loc.cit.}]). 

Then the {\em intersection}
 $\xi\mapsto \xi\cap Y_{\alpha}$ is defined as follows; 
$\xi$ {\em determines an element in
 $H_k(|\xi|, |\partial\xi|)$ via Fact, and the image of that element
by the above composition determines $\xi\cap Y_{\alpha}$, }
again by the  above Fact. 

\noindent {\bf Remark.}\ It is important to note that 
$ |\xi|\cap Y_{\alpha}$ is  not necessarily equal to 
$|\xi\cap Y_{\alpha}| $ since  the chains are not necessarily transversal
in the differential sense (cf. II.3.7).

\vspace{2mm}

Next, we define a (boundary) operator $
\cap:C_k^\pitchfork(\tilde Y^p)\rightarrow
C_{k-2}^\pitchfork(\tilde Y^{p+1}).$
If $\xi\in C_k^\pitchfork(\tilde Y^p)$, we write 
$\xi_{\alpha_1<\cdots<\alpha_p}^k$ for the corresponding components of 
$\xi$ corresponding to the decomposition of $\tilde{Y}^p$ (the index $k$ 
emphasizes the dimension, and sometimes it is omitted). Then we define:
$$\cap(\xi_{\alpha_1<\cdots<\alpha_p}^k)=\sum_{
\alpha_1<\cdots<\alpha_i<\alpha<\cdots<\alpha_p}
(-1)^{k+i}\ \xi_{\alpha_1<\cdots<\alpha_p}^k\cap Y_\alpha.
$$

\subsubsection{ }
\begin{lemma}
\ a) \  $\cap^2=0$\ and \ b)\  $\partial\cap+\cap \partial =0$.
\end {lemma}

\begin {proof}
a)  \ 
Let $\alpha<\beta$ and $\alpha,\beta\not\in \{\alpha_1, \ldots, \alpha_p\}$,
where $\alpha_1<\cdots<\alpha_i<\alpha<\alpha_{i+1}
<\cdots<\alpha_{i+j}<\beta<\cdots<\alpha_p$. 
>From the properties of the intersection (see, e.g. [GM], page 144),
one obtains that 
$(\xi_{\alpha_1<\cdots<\alpha_p}\cap Y_{\alpha})\cap Y_{\beta}=
(\xi_{\alpha_1<\cdots<\alpha_p}\cap Y_{\beta})\cap Y_{\alpha}.$
Then the sum of the following two contributions in 
$\cap^2$ is zero:
\[\cap_{Y_\beta}(\cap_{Y_\alpha}
\xi_{\alpha_1<\ldots<\alpha_p}^k)=(-1)^{k+i
}\cap_{Y_\beta}(\xi_{\alpha_1<\cdots<\alpha_p}^k
\cap Y_\alpha)=(-1)^{j+1}(\xi_{\alpha_1<\cdots<\alpha_p}^k\cap
Y_\alpha)\cap Y_\beta\]
\[\cap_{Y_\alpha}(\cap_{Y_\beta}
\xi_{\alpha_1<\ldots<\alpha_p}^k)=(-1)^{k+i+j}
\cap_{Y_\alpha}(\xi_{\alpha_1 <\cdots
<\alpha_{i+j}<\cdots <\alpha_p}^k\cap
Y_\beta)=(-1)^{j}(\xi_{\alpha_1<\cdots<\alpha_p}^k\cap
Y_\beta)\cap Y_\alpha.\]

b)  follows from the  definition of $\cap$ and $\partial$, 
including the sign $(-1)^k$ when $ \dim |\xi|= k$.
\end {proof}

\vspace{2mm}

The fact that 
the homology of $(C_{*}^{\pitchfork}(\tilde Y^p),\partial)$ is
exactly $H_{*}(\tilde Y^p)$ (cf. Lemma II.2.1) and part (b) of II.2.4 
show that
$\cap$ induces an operator $H_{*}(\tilde Y^p)\rightarrow
H_{{*}-2}(\tilde Y^{p+1})$. 

\subsubsection { }
\begin{lemma}
For any $p\geq 0$, the operator $H_{*}(\tilde Y^p)\rightarrow
H_{{*}-2}(\tilde Y^{p+1})$  induced by  $\cap$ is 
$$
\cap[\xi_{\alpha_1<\cdots<\alpha_p}^k]=\sum_{
\alpha_1<\cdots<\alpha_i<\alpha<\cdots<\alpha_p}
(-1)^{k+i}\ [\xi_{\alpha_1<\cdots<\alpha_p}^k]\cap[Y_\alpha]
$$
where 
$\cap[Y_\alpha]$ denotes the homological Gysin map, or transfer
map $i_{!}$, where 
$i:Y_{\alpha_1,\ldots,\alpha_p}\cap Y_{\alpha}\hookrightarrow
Y_{\alpha_1,\ldots,\alpha_p}$ is the natural inclusion.
More precisely, 
$i_!=PD\circ i^*\circ PD$, where $PD$ denotes the Poincar\'e Dualities
 in the  corresponding spaces (cf. e.g.  [Br], page 368). 
\end{lemma}

\begin{proof}
The result follows from the definition of the intersection 
$\xi\cap Y_{\alpha}$.
\end{proof}

\subsubsection{}{\bf The Poincar\'e duality map.} \ 
We will need later the
Poincar\'e isomorphism between homology of $ Y $ and cohomology of $ X $
with support in $ Y $, so we adopt here to our situation the
construction  of  [GM], pages 139-140, inspired by  the classical construction
which provides the Poincar\'e Duality for manifolds (see e.g. [Br],
page 338).

Similarly as above, for any triangulation $T$, compatible with the 
stratification, one can consider the chain complex of simplicial
cochains $(C^*_T(\tilde{Y}^p),\delta)$ of $\tilde{Y}^p$.
Here  $C^i_T(\tilde{Y}^p)=Hom(C_i^T(\tilde{Y}^p),\Z)$.
Let $T'$ be the first barycentric subdivision of $T$, and let $\hat{\sigma}$
 denote the barycentre of the simplex $\sigma\in T$. Let $T_i$ be the 
$i$-skeleton of $T$, thought of as a subcomplex of $T'$. It is spanned by all
vertices of   $\hat{\sigma}$   such that $\dim \sigma\leq i$.
Let $D_i$ be the $i$-coskeleton spanned, as a subcomplex of $T'$, 
by all the 
vertices  $\hat{\sigma}$ such that $\dim\sigma\geq i$.
There are canonical simplex preserving deformation retracts:
\subsubsection{}\hspace{2cm}
$X-|T_i|\to |D_{i+1}|$\ \ \and \ $X-|D_{i+1}|\to |T_i|$.

\vspace{2mm}

Now, identify $C^i_T(X)$ with $\oplus_{\dim\sigma=i}H^i(\sigma,\partial\sigma)
=H^i(|T_i|,|T_{i-1}|)$, and define
$pd:C^i_T(X)\to C^{T'}_{2n-i}(X)$ by the following composition:
$$H^i(|T_i|,|T_{i-1}|)$$
$$\ \ \ \ \ \downarrow\ \cap[X]$$
$$H_{2n-i}(X-|T_{i-1}|,X-|T_i|)$$
$$\hspace{3.5cm} \approx\ \downarrow\ \mbox{(deformation retract)}$$
$$H_{2n-i}(|D_i|,|D_{i+1}|)$$
$$\downarrow$$
$$H_{2n-i}(|T'_{2n-i}|,|T'_{2n-i-1}|)$$

Actually, $H_{2n-i}(|T'_{2n-i}|,|T'_{2n-i-1}|)=C^{T'}_{2n-i}(X)$, but we have
even  something more.  Since the image of any chain by this composition is
supported by union of $|D_i|$'s, and any $|D_i|$ is dimensionally transverse,
one obtains a homomorphism $pd:C^i_T(X)\to C_{2n-i}^{T',\p}(X)$. 
Similarly, for any $\Y$, one can define a homomorphism:
$$pd:C^i_T(\Y)\to C_{2n-2p-i}^{T',\p}(\Y).$$
By [GM] (7.2), this is a chain map:

\subsubsection{}\  \hspace{2cm} $\partial\circ pd=pd\circ \delta$.

\subsubsection{}{\bf Lemma.} \ {\em a)\ 
 \ Fix $\alpha\not\in \{\alpha_1,\ldots,
\alpha_p\}$. Then the following diagram is commutative:}

$$\begin{array}{ccc}
C^i_T(\Y)&\stackrel{pd}{\longrightarrow}&
C_{2n-2p-i}^{T',\p}(\Y)\\
\Big\downarrow\vcenter{\rlap{$i^*$}}&&
\Big\downarrow\vcenter{\rlap{$\cap$}}\\
C^i_T(\Y\cap Y_{\alpha})&\stackrel{pd}{\longrightarrow}&
C_{2n-2p-i-2}^{T',\p}(\Y\cap Y_{\alpha})\end{array}$$

b)\ {\em  The above diagram (and II.2.8) provides at homology level
a  commutative diagram: 

$$\begin{array}{ccc}
H^i(\Y)&\stackrel{PD}{\longrightarrow}&
H_{2n-2p-i}(\Y)\\
\Big\downarrow\vcenter{\rlap{$i^*$}}&&
\Big\downarrow\vcenter{\rlap{$\cap$}}\\
H^i(\Y\cap Y_{\alpha})&\stackrel{PD}{\longrightarrow}&
H_{2n-2p-i-2}(\Y\cap Y_{\alpha})\end{array}$$
where the horizontal maps are the Poincar\'e Duality isomorphisms
(cf. also with II.2.5).}
\begin{proof} Use the definition of \  $\cap$ and $pd$ and 
II.2.7.\end{proof}

\subsection{}{\bf The chain correspondence $\Pi^{-1}_{\star}$.} 

\subsubsection{} {\bf Triangulation of $ Z$ }. 
Let us indicate quickly  that the real blow up $ \Pi: Z \to X$ of
$ Y$ in $ X$ can be realized as a semi-analytic subset of a real vector
space.
 Since $ X$ is compact, there is
no difference between bounded semi-analytic and
subanalytic simplices  [H2, Def. 3.3, page  179].
 The variety $ X$ is locally, analytically,
isomorphic to a ball in
$\C ^n$, such that  the variety $ Y_J$ is defined by $   
z_1\cdot\ldots\cdot z_j = 0 $ where
$j = \vert J \vert$,  then the variety $ Z$ is
isomorphic to the    semi-analytic subset of a real vector space: $
T^n \times (\R^+)^n $ with the
the product structure
 of $ (\R^+)^n $ and the $n$-dimensional  torus. The projection $ \Pi
$ is written, in the coordinates $ \theta_k$ on the torus,
$\rho _k$ on the reals  and $ z_k $  on $
\C^n$ for $ k = 1,...,n $,  as $ z_k \circ \Pi = \rho _k \
e^{i\theta _k} $.  We have
$$ \Pi^{-1} ( Y_J - Y^{j+1} ) \simeq  T^j\times
T^{n-j}\times (\R^*)^{n-j} \simeq T^j\times (\C^*)^{n-j}$$
 where we use in the last  isomorphism the correspondence
defined by $ z_k = \rho _k e^{i\theta _k}, k = j+1,...,n$ to
describe $ Y_J - Y^{j+1}$.

Define the strict transform
of a  simplex
$\Delta$ on $ \tilde Y^r $ as the closure of the inverse image 
in  $ Z$  of the  
complement of $Y^{r+1} $ in $\Delta$.

Let $\Delta _\alpha$ be a semi-analytic (resp. sub-analytic) simplex
contained in the domain
of a chart  on $ X$ and defined by a set of equalities and
inequalities of analytic functions
$f_{\alpha,\lambda} (x.,y.)_{\lambda \in \Lambda}$ in the
real coordinates $(x.,y.)$ of $ \C^n$ ; its inverse image is
obtained by replacing $ ( x.,y. )$ with $ \rho _. ( cos (
\theta .),  sin (\theta .)) $ in the defining equations.
 The strict transform of the simplices $\Delta_\alpha$
in $ X$  are
semi-analytic  subsets of $ Z $ (resp. sub-analytic).

They  define  in $ Z $ a decomposition in a finite system of
sub-analytic sets, so that we can use Hironaka's existence
theorem [H2, page 180] to define  a triangulation by sub-analytic
simplices, subordinated to this decomposition. 

\subsubsection{}
 In this  subsection we present our {\em
main topological tool}: a  correspondence
denoted $\Pi^{-1}_{\star}$
 which gives   closed chains in $ Z $, respectively in its
boundary
$\partial Z$, as ``strict inverse''   of chains on
$X$. But first we have to
fix some orientation conventions
regarding the fibers $\Pi^{-1}(y^o)$ for different points  $y^o\in Y$. 

We denote the oriented boundary  of a disc in the complex plane
by $S^1=\partial D$  (where we consider the natural orientation as a boundary).
For any point  $y^o\in Y^1-Y^2$, the circle  $\Pi^{-1}(y^o)$
appears as the boundary of the complement of a disc in  $\C$, hence  
it is $-S^1$.     If $y^o\in Y^p-Y^{p+1}$, the situation is similar.  
In a local model,
    when $Y $ is defined by $ \{y_1\ldots y_p=0\}\subset\C^p$ where $y^o$ 
stays for the 
origin. Let $U $ be defined by $\{y\, |\, \min |y_i|\leq 1\} \subset\C^p$
    and the component $Y_{\alpha_i}$ by $y_i=0$. Then
    $\Pi^{-1}(y^o)$, set-theoretically, is the tori $S_1^1\times\ldots\times
    S_p^1$,  where $S_i^1 $ is defined by $\{y_i:\left|y_i\right|=1\}$.  
    Then, we fix the 
    orientation of $\Pi^{-1}(y^o)$ as given by the product orientation
    $(-S_1^1)\times\cdots\times (-S_p^1)$.
Moreover, if $\sigma$ is a contractible $C^{\infty}$-submanifold in 
    $Y_{\alpha_1}\cap\cdots\cap Y_{\alpha_p}-Y^{p+1}$
    $(\alpha_1<\cdots<\alpha_p)$, then on the set-theoretical
    inverse image  $ \Pi^{-1}(\sigma)$  we define the product orientation
    $\sigma\times (-S_{\alpha_1}^1)\times\cdots\times
    (-S_{\alpha_p}^1)$.

With this notations, the following holds.

\subsubsection{}
\begin{lemma}
Fix indices
$\alpha_1<\cdots<\alpha_i<\alpha<\alpha_{i+1}<\cdots<\alpha_p$ and let
$\sigma\subset Y_{\alpha_1}\cap\cdots\cap Y_{\alpha_p}$ be an oriented 
${\C}^\infty$ sub-manifold with boundaries such that $\sigma$ and
$\partial\sigma$ intersect $Y_\alpha$ transversally, but $\sigma\cap
Y_{\beta}=\emptyset$ for any $\beta\not\in 
\{\alpha,\alpha_1,\ldots,\alpha_p\}$.   
(The $C^{\infty}$ transversal intersection is considered in
$Y_{\alpha_1}\cap\cdots \cap Y_{\alpha_p}$.) Then
\[\partial(\Pi^{-1}\sigma)=
\Pi^{-1}(\partial\sigma) +(-1)^{\dim\sigma+i}\Pi^{-1}(\sigma\cap Y_{\alpha}).\]
\end{lemma}

\begin{proof}

In a neighbourhood of an intersection point $p\in(\sigma\cap Y_\alpha)$
($\resp(\partial\sigma)\cap Y_\alpha),\ \  \sigma$ has a product
structure $D\times T$ where $T$ is a ball ($\resp$ a half ball) in
$Y_\alpha \cap Y_{\alpha_1} \cap \cdots \cap Y_{\alpha_p}$, and $D$ is a
real 2-disc transversal to $Y_\alpha$.  We denote the boundary
of $D$ by $\partial D$. Let $D_\eta =D-\{$small open disc
of radius  $\eta\ $  and origin $\,0\}$;  i.e.  $\partial D_{\eta}=
\partial D -S_{\eta}^1$.  Then we have the diffeomorphism
$\Pi^{-1}(D\times T)\approx  \Pi^{-1}(D_\eta\times T)$. 
Since $ \Pi^{-1}(D_\eta\times T) = D_\eta\times T\times \Omega$, 
 where $ \Omega  = (-1)^{p}S_{\alpha_1}^1\times\cdots
\times S_{\alpha_p}^1 $, we get: $\partial\Pi^{-1}(D\times T)
= \partial \Pi^{-1}(D_\eta\times T)
=\partial D\times T\times \Omega   -S_\eta^1\times
T\times \Omega +  D_\eta\times\partial T\times \Omega=$

$=(-1)^{\dim T+i}
T\times(-1)^{p+1}S_{\alpha_1}^1\times\cdots\times
S_\alpha^1\times\cdots\times S_{\alpha_p}^1+ \partial (D\times T)
 \times \Omega $

$= (-1)^{i+\dim\sigma}\Pi^{-1}(\sigma\cap Y_\alpha)+\Pi^{-1}\partial\sigma.$
\end{proof}

\subsubsection{} Next, we want to lift an arbitrary sub-analytic  chain.
First notice that  for each  sub-analytic subset  $S\subset \tilde Y^r$, 
the strict transform $\tilde{S}:= cl ( \Pi^{-1} (S - Y^{r+1} ))$ (i.e.
the closure of the inverse image of the
complement of $Y^{r+1} $ in $S$) is  a sub-analytic subset of $Z$. 
Indeed, the strict transform can be constructed using a real
analytic isomorphism and  the permitted
operations listed in II.1.2 (2-3) (cf. also with II.3.1).
 
Fix again the integer  $p\geq 0$,  and for any chain
$\xi\in C^{\pitchfork}_k(\tilde{Y}^p)$ we want to construct a chain
$\Pi^{-1}_p(\xi)$ in $C_{k+p}(\partial Z)$ if $p\geq 1$, respectively in
$C_k(Z)$ if $p=0$.   By the above remark, 
the inverse image $\Pi^{-1}(|\xi|)$ of the support of $\xi$ 
is a $(k+p)$--dimensional  sub--analytic  subset of $ Z $.
Similarly as above, in the construction we will use Fact II.2.3.

\vspace{2mm}

Like in the case of Leray's  cohomological residue, we have the following
morphism of relative homology:
$$ \Pi_{p,rel}^{-1}:
H_k(|\xi|,|\partial\xi|\cup (|\xi|\cap Y^{p+1})) \longrightarrow 
 H_{k+p}(\Pi^{-1}(|\xi|),
\Pi^{-1}(|\partial \xi|)\cup \Pi^{-1}(|\xi|\cap Y^{p+1})).$$
Since $\Pi^{-1}(|\xi|)\to
|\xi|$ is an oriented  fiber  bundle over $|\xi|-Y^{p+1}$,
with fibers $(S^1)^p$, the classical spectral sequence argument leads to an  
isomorphism. Indeed, we use deformation retract tubular neighbourhoods to 
thicken $ |\xi|\cap Y^{p+1} $ and its inverse image 
$\Pi^{-1}(|\xi|\cap Y^{p+1})$, then the excision theorem to reduce 
to the fiber bundle case and show that the relative morphism above is 
an isomorphism. 

Next,  we use the morphism
$$ i :  H_k(|\xi|,|\partial\xi|)
\longrightarrow H_k(|\xi|,|\partial\xi|\cup (|\xi|\cap Y^{p+1}))$$ 
induced by the inclusion
$(|\xi|,|\partial\xi|)\longrightarrow (|\xi|,|\partial\xi|\cup(|\xi|
\cap Y^{p+1}))$. Notice that
$i$ is an isomorphism since codim$|\xi|\cap Y^{p+1}$ is $ 2 $ in $ |\xi|$ and
in the long homology exact sequence of the pair 
$H_i(|\partial \xi|\cup (|\xi|\cap Y^{p+1}),|\partial \xi|)
= H_i(|\xi|\cap Y^{p+1},|\partial \xi|\cap Y^{p+1})=0$ if $i=k$ or
$k-1$.

Now, the image of $\xi$, via the composed map $\Pi_{p,rel}^{-1}\circ i$,
 determines completely $\Pi^{-1}_p(\xi)$ (via Fact). By definition, 
this is the  application $\xi\mapsto \Pi^{-1}_p(\xi)$.  
Sometimes the index $p$ will be replaced by  $\s$.

\vspace{2mm} 

In a simple way, this chain $\Pi^{-1}_p(\xi)$ can be defined  as follows. 
Write $\xi$ as  a finite sum $\sum a_{\sigma} \,\sigma$, where $\sigma$ are 
the simplices  in the support of $\xi$  with non-zero coefficients.
 Let $\sigma^0$ be
$\sigma-Y^{p+1}$. Then the closure $cl(\Pi^{-1}(\sigma^0))$ of
$\Pi^{-1}(\sigma^0)$ is a $(k+p)$--dimensional sub--analytic set.
Then define $\Pi^{-1}_p(\xi)$ by $\sum\ a_{\sigma}\,cl(\Pi^{-1}(\sigma^0))$.

\vspace{2mm}

The next result generalizes the above lemma.

\subsubsection{} 
\begin{prop}
Fix the integer $p\geq 0$ and indices 
$\alpha_1<\cdots<\alpha_p$. For any $\alpha\not\in \{
\alpha_1,\ldots,\alpha_p\}$, set $i(\alpha)=i$ if 
$\alpha_1<\cdots<\alpha_i<\alpha<\alpha_{i+1}<\cdots<\alpha_p$.
Then for any $\xi\in C^{\pitchfork}_k(Y_{\alpha_1,\ldots,\alpha_p})$:

\[\partial(\Pi^{-1}_p\xi)=
\Pi^{-1}_p(\partial\xi)+\sum_{\alpha}\ 
(-1)^{k+i(\alpha)}\ \Pi^{-1}_{p+1}(\xi\cap Y_{\alpha}).\]
Above,  the sum is over all $\alpha$ with 
$\alpha\not\in \{\alpha_1,\ldots,\alpha_p\}$.
The equality is considered in  $C_*(\partial Z)$ if $p\geq 1$, and in 
$C_*(Z)$ if $p=0$.
\end{prop}
\begin{proof}\ The contribution $\Pi^{-1}_p(\partial\xi)$ is clear. 
Next we want to determine the coefficients  of the simplices
which lie in $\Pi^{-1}(|\xi|\cap Y_{\alpha})$. 
It is enough to work modulo
$Y_{\alpha}\cap (\cup_{\beta}Y_{\beta})$, where the union is over 
all $\beta\not\in \{\alpha,\alpha_1,\ldots,\alpha_p\}$, 
since it intersects $ |\xi| $ in codimension 4 and its
inverse image intersects $ \Pi^{-1}(|\xi|) $ in codimension 2. 

Consider the composition:
$$H_{k+p}(\Pi^{-1}(|\xi|),
\Pi^{-1}(|\partial \xi|)\cup \Pi^{-1}(|\xi|\cap Y^{p+1}))$$
$$\hspace{.1cm}\downarrow\ \partial$$
$$H_{k+p-1}(\Pi^{-1}(|\partial \xi|\cup(|\xi|\cap Y^{p+1})),
\Pi^{-1}(|\partial \xi|\cup(|\xi|\cap \cup_{\beta}Y_{\beta})))$$
$$\approx\ \downarrow\ e\ \ $$
$$H_{k+p-1}(\Pi^{-1}(|\xi|\cap Y_{\alpha}),
\Pi^{-1}(|\partial \xi|\cap Y_{\alpha})\cup \Pi^{-1}
(|\xi|\cap Y_{\alpha}\cap\cup_{\beta}Y_{\beta}))).$$
Again, the union $\cup_{\beta}$ is over all $\beta\not\in
\{\alpha,\alpha_1,\ldots,\alpha_p\}$.  If $A=\Pi^{-1}(|\xi|)$,
$B=
\Pi^{-1}(|\partial \xi|\cup(|\xi|\cap \cup_{\beta}Y_{\beta})))$, 
and $C=
\Pi^{-1}(|\xi|\cap Y_{\alpha})$, then the first map is the boundary operator
$H_{k+p}(A,B\cup C)\to H_{k+p-1}(B\cup C,B)$, and the second is the
excision isomorphism $H_{k+p-1}(B\cup C,B)\to 
H_{k+p-1}(C,B\cap C)$. 

The composed map  $e\circ \partial$ and the map $\xi\mapsto \xi\cap Y_{\alpha}$
are connected by the following diagram, which commutes up to sign:

$$H_k(|\xi|,|\partial \xi|)
\hspace{.5cm}\stackrel{\Pi^{-1}_p}{\longrightarrow}\hspace{.5cm}
H_{k+p}(\Pi^{-1}(|\xi|),
\Pi^{-1}((|\partial \xi|)\cup (|\xi|\cap Y^{p+1})))$$
$$\hspace{.1cm} \downarrow \ \cap Y_{\alpha}\hspace{6cm}
\downarrow \ e\circ \partial\hspace{2cm}$$
$$H_{k-2}(|\xi|\cap Y_{\alpha},|\partial \xi| \cap Y_{\alpha})
\stackrel{\Pi^{-1}_{p+1}}{\longrightarrow}
H_{k+p-1}(\Pi^{-1}(|\xi|\cap Y_{\alpha}),
\Pi^{-1}((|\partial \xi|\cap Y_{\alpha})\cup
(|\xi|\cap Y_{\alpha}\cap\cup_{\beta}Y_{\beta})))$$

The commutativity of the diagram (up to a sign), basically comes from the fact
that for a manifold with boundary, 
the Lefschetz duality identifies the boundary operator
in homology with the restriction map (to the boundary) in cohomology
(see e.g. [Br], page 357). 
The above sign is universal, depends only on the orientation conventions. 
Therefore, it can be determined using $C^{\infty}$-transversal chains, as in 
lemma  II.3.3.\end{proof}

The above proposition together with the definition of the operator 
$\cap:C^{\pitchfork}_k(\tilde{Y}^p)\to
C^{\pitchfork}_{k-2}(\tilde{Y}^{p+1})$, have the following corollary.

\subsubsection{} {\bf Corollary.}  \ 
$\partial \Pi^{-1}_\s=\Pi^{-1}_\s(\partial+\cap).$

\subsubsection{}{\bf Example.}\  Let $\xi\in  C_*^{\pitchfork}(X)$. 
In general, the three sets $\Pi^{-1}(|\xi|\cap Y)$ (here $\Pi^{-1}$ denotes
the set-theoretical inverse image), $|\Pi^{-1}_\s(\xi)|\cap \partial Z$ and
$|\partial \Pi^{-1}_\s(\xi)|$ are all distinct. To see this,
set $X=\C^2$, $Y=\{0\}\times \C$ ($Y^2=\emptyset$). 
Fix coordinates $z_j=x_j+iy_j$
($j=1,2$ ) in $\C^2$. Let $\xi$ be a chain with support 
$(x_1-1)^2+y_1^2+x_2^2-1=y_2=0$, and coefficient one. 

Then $|\xi|$  is a 2-dimensional real sphere with $|\xi|\cap Y=(0,0)$, hence 
$\Pi^{-1}(|\xi|\cap Y)=S^1$. On the other hand, the intersection chain 
$\xi\cap Y=0$  (even if $|\xi|\cap Y\not=\emptyset$). Actually, $\xi\cap Y$ 
is the chain supported by $|\xi|\cap Y$  with coefficient the intersection
multiplicity of $|\xi|$ and $Y$, which is zero. Therefore, $\partial 
\Pi^{-1}_\s
(\xi)=\Pi^{-1}_\s\partial \xi=0$, hence $|\partial \Pi^{-1}_\s(\xi)|
=\emptyset$. 
Finally, the reader is invited to verify that $|\Pi^{-1}_\s(\xi)|\cap 
\partial Z$ is a {\em half circle} in $S^1=\Pi^{-1}(0,0)$. 

\subsubsection{}{\bf Remark.}\ In the above discussion, 
corresponding to the integer $p=0$, one can replace  the space
$X$ by the (compact) tubular neighbourhood $U=U_{\alpha}=\Pi(Z_{\alpha})$ (cf.
I.1.2). This means that the complex $C_*(X)$ is replaced by $C_*(U)$,
and $C^{\p}_*(X)$ by $C_*^{\p}(U)$.  Notice that $\cap: C^{\p}_*(U)\to
C^{\p}_{*-2}(\tilde{Y}^1)$  is well-defined and 
still satisfies $\cap\partial+\partial\cap=0$.
Moreover, one has a map $\Pi^{-1}_\s:C^{\p}_*(U)\to C_{*}(Z_{\alpha})$, 
defined similarly as 
$\Pi^{-1}_\s:C^{\p}_*(X)\to C_{*}(Z)$, which satisfies 
$\partial \Pi^{-1}_\s=\Pi^{-1}_{\s}(\partial+\cap)$.

\subsection{The complex of sheaves.}

\noindent   $C_*(\tilde{Y}^p)$ has a natural sheafification 
$\c{C}^{*}_{\tilde{Y}^p}$ which is fine  ([GM2] page 97, or [B] page 33). 
Actually, the complex 
$\c{C}^{*}_{\tilde{Y}^p}$ is quasi-isomorphic to the dualizing complex 
$\c{D}^*_{\tilde{Y}^p}$ (cf. [GM2] page 97, or [B] page 33).

What we will need in this paper, is a different construction, we will 
sfeafify the dual complex, and we will obtain a quasi-isomorphism with 
$\Z_{\tilde{Y}^p}$. 

For any open set $V\subset \tilde{Y}^p$, let $C_k(V,\tilde{Y}^p)$ be the 
subgroup of chains $\xi\in C_k(\y)$ with compact support $|\xi|$ in $V$.
Obviously, for any open pair $V\subset W\subset \y$, there is a natural
inclusion $C_k(V,\y)\to C_k(W,\y)$. Now, define the dual 
$C^k(V,\y)$ by $Hom_{\Z}(C_k(V,\y),\Z)$. Then for any $V\subset W$ as above,
the ``restriction'' $C^k(W,\y)\to C^k(V,\y)$ defines a presheaf $\c{C}^k(\y)$
 on $\y$, satisfying the condition (S2) (i.e. it is ``conjuctive'' in the 
terminology of [Br]). Let $\bar{\c{C}}^k(\y)$ be the associated sheaf with 
global sections $\bar{C}^k(\y)$; and let $C^k_0(\y)$ be the subgroups of
elements of $C^k(\y)$ with empty support. Then 
$$0\to C_0^k(\y)\to C^k(\y)\to \bar{C}^k(\y)\to 0$$
is exact (cf. [Br] page 22). Moreover, $C^*(\y)$ and $C_0^*(\y)$ form
 complexes,
and $H^*(C^*_0(\y))=0$ (by subdivision argument, see [Br2], page 26
in the case of singular chains). 

Therefore, the complexes $C^*(\y)$ and $\bar{C}^*(\y)$ are quasi-isomorphic.
On the other hand, for any open $V$, there is an augmentation map 
$C_0(V,\y)\to \Z$ which give rise to a resolution $0\to
\Z_{\y}\to \bar{\c{C}}^*(\y)$. The sheaf $\bar{\c{C}}^k(\y)$ is a module over
the ring of $\Z$-constructible  functions on $\y$.
Indeed,  for any  constructible function $f$ and $\varphi\in C^k(\y)$
one can define $f\cdot \varphi\in C^k(\y)$ by
$(f\cdot\varphi)(\sigma)=f(\hat{\sigma})\varphi(\sigma)$,
where $\hat{\sigma}$ is the barycenter of  $\sigma$. 
This shows that the above resolution is a resolution of fine sheaves.

\section{Double complexes and their spectral sequences}

We construct in this section
various homological complexes giving rise to  various homology groups
with their weight filtration. The cases $(Y,\emptyset)$ and 
$(X,Y)$ are well-known, but we will need them in the construction of the
double complex of $\partial U$, so we include them in our presentation
as well. The new result is the construction of the complexes of 
the pairs $(X-Y,\emptyset)$, $(X,X-Y)$ and $(\partial U,\emptyset)$,
in which cases we will use the dimensionally transverse cycles. 
The most involved  case is $\partial U$, which is separated in section VI.
In the next paragraphs, as an introductory guide, 
we will stress the case $X-Y$.

The {\em main new object is homological double complex} of $X-Y$:
$$A_{s,t}(X-Y):=C_{t+2s}^\pitchfork (\tilde Y^{-s}), s\leq 0, t+2s\geq 0$$
with 
 $D=\partial + \cap$ as the  differential of the total complex $\Tot_*
(A_{**}(X-Y))$.

Our aim, reformulated for the case of $X-Y$,
 is to prove that the homology of the above total complex is
the homology of $ X - Y$, and its weight filtration provides Deligne's
weight filtration on $H_*(X-Y,\Q)$. More precisely:

\vspace{1mm}

\noindent III.0. {\bf Theorem.}\ {\em  
Let $Y=\bigcup_{i\in I} Y_i$ be a NCD in a smooth proper algebraic variety
$X$ over $\C$ and $\Pi: Z \rightarrow  X $ the  projection
 as above. The generalized Leray inverse image defines a
quasi-isomorphism
$$ \Pi^{-1}_\s: \Tot_* (A_{**}(X-Y)) \longrightarrow C_*(Z).$$
In other words, a collection of 
dimensionally transversal  geometric chains $\{c_r\}_{r\geq p}$, where $c_r$
is a chain on  $\tilde Y^r$ satisfying $\partial c_p=0$ and
$\partial c_{r+1}=c_r\cap\tilde Y^{r+1}$ for any $r\geq p$ ($p\geq 0$), 
defines a cycle $c=\sum_{r\geq p} \Pi_r^{-1}c_r$ in
$ Z $. Moreover, their homology classes generate $H_*(Z,\Q)$. 

Geometrically this  means the following: the highest nontrivial weight in 
$H_k(X-Y,\Q)$ is $-k$, and for any $p\geq 0$,
$W_{-k-p}H_k(X-Y,\Q)$ is generated by classes of cycles $e$, with a 
chain decomposition $e=\sum_{r\geq p}e_r$ such that
each $e_r$ is the closure of  a ``fiber bundle'' over a $c_r-Y^{r+1}$ 
with fiber $(S^1)^r$, where $c_r$ are chains  in $\tilde{Y}^r$ and 
$c_r-Y^{r+1}$ is considered in $Y^r$. 
This is realized geometrically in a small tubular
neighbourhood of $Y^r$. 
The fact that the collection $\{c_r\}_r$ satisfies the above 
condition is equivalent with the fact that $\{e_r\}_r$ fit together without 
boundary, forming a closed cycle. 
}

\vspace{1mm}

Above we identified $H_*(C_*(Z))$ and $H_*(X-Y)$. For more details, see III.4
and V.2. 

\vspace{2mm}

 This result can be obtained by duality  with Deligne's
 logarithmic differential complex.  Here  we  follow a purely topological
path:

\vspace{1mm}

1) We start with the basic Mayer--Vietoris resolution  for the 
$NC$-divisor $Y$.

2) We introduce the  relative homology of $ X \, mod \, Y$ using the 
classical cone construction, called here {\em mixed cone}, in order to carry
the  construction with weight filtrations (a diagonalisation process for
the weights is needed  in order to stay in the category of mixed Hodge
complexes [D] so to get mixed Hodge structures on their  (co)homology by
the basic result of Deligne). 

3)  We give a topological construction of Poincar\'e duality  between
the cohomology of $ Y$ and the homology of the pair $ ( X, X - Y ) $.

4) The above duality extends to a duality  between the cohomology of
the pair $ ( X , Y ) $   and the homology of $ X - Y $. This provides
 the wanted result.

Moreover we obtain the
degeneration of the weight spectral sequence on the total complex of
dimensionally transversal chains.

\vspace{2mm}

\noindent {\bf Notations and preliminary remarks.} 
For any double complex $(A_{**},\partial,\delta)$, 
we denote its total complex by $(Tot_*(A_{**}), D)$, where $Tot_k(A_{**})=
\oplus_{s+t=k}A_{st}$ and $D=\partial + \delta$. (Here the degree of
$\partial$ is $(0,-1)$,  of  $\delta $ is  $(-1,0)$.)
The weight filtration of $A_{**}$ is defined by 
$W(A_{**})_s:=\oplus_{p\leq s}A_{pq}$. 
The homological spectral sequence associated with the weight filtration $W$
is denoted by $E^r_{**}$.
Recall that $E^1_{st}=H_{t}(A_{s*},\partial)$ and $d^1$ is induced by 
$\delta$. However, we will violate the weight notation on the $\infty$-term,
and we will use Deligne's convention: on $E_{st}^{\infty}$ the weight 
is $-t$ (instead of $s$); i.e. image$\{H_i(Tot_*(W_s)\to H_i(Tot_*(A))\}$ 
is $W_{s-i}H_i(Tot_*(A))$.

The {\em dual double complex} of $A_{**}$ is $B_{st}=
Hom_{\Z}(A_{st},\Z)$ with the corresponding dual maps. Its weight 
filtration is
$W(B)_{s}=\{\varphi\in A^*:\varphi(W(A)_{-s-1})=0\}$, or
equivalently, $W(B)_{-s}=\oplus_{p\geq s}B_{pq}$. 

For the definition and properties of  ``cohomological mixed Hodge complexes''
see [D.III]. In the next proofs the notion of ``mixed cone'' will be
important, this corresponds to the mapping cone in the category of 
mixed Hodge complexes. For details, see [D. III], page 21, or [E], page 49.

\subsection{}{\bf The homological double complex of $Y$ 
(the homology of the $NCD$)}.\\

Consider the double complex $A_{s,t}(Y):=C_t(\tilde{Y}^{s+1})$ 
(with $s\geq 0$ and
$t\geq 0$) together with the natural operators:
$$\partial: C_k(\tilde{Y}^p)\to C_{k-1}(\tilde{Y}^p)\ \mbox{and}
\ \ i: C_k(\tilde{Y}^p)\to C_{k}(\tilde{Y}^{p-1}).$$
Here $i$ is defined as follows. If $\oplus_{\alpha_1<\cdots<\alpha_p}
c_{\alpha_1,\ldots,\alpha_p}\in 
C_k(\tilde{Y}^p)$, then $i(\oplus_{\alpha}c_{\alpha})=\oplus_{\beta}d_{\beta}$
if:
$$d_{\alpha_1,\ldots,\alpha_{p-1}}=\sum_{\alpha_1<\cdots<\alpha_i<\alpha<\cdots
<\alpha_{p-1}}\ (-1)^{k+i}\ 
i_{\alpha_1,\ldots,\alpha_{p-1};\alpha}(c_{\alpha_1,\ldots,\alpha_i,\alpha,
\ldots,\alpha_{p-1}}),$$
where $i_{\alpha_1,\ldots,\alpha_{p-1};\alpha}$
is the natural inclusion 
$Y_{\alpha_1,\ldots,\alpha_i,\alpha,
\ldots,\alpha_{p-1}}\hookrightarrow
Y_{\alpha_1,\ldots,\alpha_{p-1}}.$

\subsubsection{}{\bf Lemma.} \ {\em a)\ $i^2=0$, and b)\ $i\partial+\partial i=0$.}

\vspace{2mm}

In particular, $D:=i+\partial$  is a differential of
the total complex $Tot_*(A_{**}(Y))$. 
The weight filtration $\{W(A)_s\}_s$ provides a spectral sequence over $\Z$.
Here are some of its properties.

\subsubsection{}{\bf Proposition.} \ {\em 

a)\ $E^1_{st}=H_{t}(\tilde{Y}^{s+1})$ and $d^1:E^1_{st}=H_t(\tilde{Y}^{s+1})
\to H_t(\tilde{Y}^s)=E^1_{s-1,t}$ is $i_*$ induced by $i$.

b)\ $E^{r}_{st}\Longrightarrow H_{s+t}(Y,\Z) $, inducing a weight filtration
on $H_*(Y,\Z)$. Actually, 
$E^{\infty}_{st}\otimes \Q=Gr^W_{-t}H_{s+t}(Y,\Q)$
(in Deligne's weight notation). 

c)\ $ E^*_{**}\otimes \Q$ 
degenerates at level two, i.e. $d^r\otimes 1_{\Q}=0$ for $r\geq 2$. }
\begin{proof}
Let $n:\tilde{Y}^p\to Y$ be the natural map. 
Recall that the $K_{\Q}^{\cdot}$ term in Deligne's cohomological mixed 
Hodge complex  associated with the space $Y$ is the ``Mayer-Vietoris
resolution'' $n_*\Q_{\tilde{Y}^{\cdot}}$: 
$$0\to n_*\Q_{\tilde{Y}^1}\to n_*\Q_{\tilde{Y}^2}\to \cdots.$$
This is considered with  its ``b\^ete'' filtration $W_{-s}
(n_*\Q_{\tilde{Y}^{\cdot}})=\sigma_{\geq s}(n_*\Q_{\tilde{Y}^{\cdot}})$.

Consider now the dual double complex  $B_{**}$ of $A_{**}$.
Then $(n_*\Q_{\tilde{Y}^{\cdot}},W)$ and 
$(B_{**},W)\otimes 1_{\Q}$  are  quasi-isomorphic. 
This follows from the discussion II.4. 
Therefore, their spectral sequence (for $r>0$) are isomorphic. 
This gives $a)$ and $b)$. 
Finally notice that by a result of Deligne, 
the weight spectral sequence (over $\Q$) 
of a  mixed Hodge complex degenerates at rank two, which provides $c)$.
\end{proof}

For a different approach and proof, see [McC], where it is proved 
that  the filtered 
dualizing complex $(i_*\c{D}_Y^*,\tau_{\leq})$ and $(i_*n_*\Q_{\tilde{Y}},
\sigma_{\geq})$ are Verdier dual (here $i:Y\to X$ stays for  the inclusion).

\subsection{}{\bf The homological double complex of $(X,Y)$.}\\

Starting with $\tilde{Y}^0=X$, we 
consider now the double complex $A_{s,t}(X,Y):=C_t(\tilde{Y}^{s})$ 
(with $s\geq 0$ and
$t\geq 0$) together with the natural operators:
$\partial: C_k(\tilde{Y}^p)\to C_{k-1}(\tilde{Y}^p)\ \mbox{and}
\ \ i: C_k(\tilde{Y}^p)\to C_{k}(\tilde{Y}^{p-1})$  as above.
Again, 
$D:=i+\partial$ is a differential of the total complex $Tot_*(A_{**}(X,Y))$. 

If we introduce 
the double complex $A_{**}(X)$ of $X$ defined by $A_{s*}(X)=C_*(X)$
if $s=0$ and $=0$ otherwise, then $A_{s*}(X,Y)=A_{s*}(X)\oplus A_{s-1,*}(Y)$, 
and $A_{**}(X,Y)$ can be interpreted  as the $Cone(\tilde{i})$ of \ 
$\tilde{i}:A_{**}(Y)\to A_{**}(X)$, where $\tilde{i}|A_{s*}(Y)=0$ if $s\not=0$,
and $\tilde{i}|A_{0*}(Y)=i:C_*(\tilde{Y}^1)\to C_*(X)$. 

\subsubsection{}{\bf Proposition.} \ {\em

a)\ $E^1_{st}=H_{t}(\tilde{Y}^{s})$ and $d^1$ is $i_*$, induced by $i$.

b)\ $E^{r}_{st}\Longrightarrow H_{s+t}(X,Y,\Z) $, inducing a weight filtration
on $H_*(X,Y,\Z)$. 
$E^{\infty}_{st}\otimes \Q=Gr^W_{-t}H_{s+t}(X,Y,\Q)$ (in 
Deligne's weight notation). 

c)\ $ E^*_{**}\otimes \Q$ 
degenerates at level two, i.e. $d^r\otimes 1_{\Q}=0$ for $r\geq 2$. }
\begin{proof}
The proof is similar as in the case of III.1.2.
In this case, $K_{\Q}^{\cdot}$ is
$$0\to n_*\Q_{\tilde{Y}^0}\to n_*\Q_{\tilde{Y}^1}\to \cdots.$$
Actually this 
(and the whole cohomologically mixed Hodge complex of $(X,Y)$) can be
constructed as a mixed cone of the complexes  of $X$, respectively of $Y$.
This is compatible with the construction of $A_{**}(X,Y)$. 
\end{proof}

\subsection{}{\bf 
The homological double complex of $(X,X-Y)$ or $(U,\partial U)$.}\\

\noindent Here $U$ is the ``tubular neighbourhood'' of $Y$. 
We define 
$$A_{s,t}(X,X-Y):=C_{t+2(s-1)}^\pitchfork (\tilde Y^{-(s-1)}),$$
with $s\leq 0$ and $t+2(s-1)\geq 0$. 
Then $\partial$ and $\cap$ act as $\partial:A_{s,t}\to A_{s,t-1}$ and
$\cap:A_{s,t}\to A_{s-1,t}$, hence $D=\partial + \cap$ is the 
differential of the total complex $\Tot_* (A_{**}(U,\partial U))$.
Corollary II.3.6 reads as:

\subsubsection{}{\bf Corollary.} {\em \ 
$\Pi^{-1}_\s:(Tot_*(A_{**}(X,X-Y)),D) \to (C_{*-1}(\partial Z),\partial)$ 
is a morphism of complexes
i.e. $\partial \Pi^{-1}_\s=\Pi^{-1}_\s D$.}

\vspace{2mm}

 Moreover, II.2.1, II.2.5 and II.2.8  and Poincar\'e duality
 imply the following result.

\subsubsection{ }\begin{prop}

a)\ $E_{s+1,t}^1=H_{t+2s}(\tilde Y^{-s})$, and 
$d_1$ is the transfer map $i_!$.

b)\ $E^{r}_{st}\Longrightarrow H_{s+t}(X,X-Y,\Z) $, and 
$E^{\infty}_{st}\otimes \Q=Gr^W_{-t}H_{s+t}(X,X-Y,\Q)$.

c)\ $ E^*_{**}\otimes \Q$ 
degenerates at level two, i.e. $d^r\otimes 1_{\Q}=0$ for $r\geq 2$. 

d)\ The Poincar\'e duality map $pd$ (cf. II.2.6)
 induces an isomorphism of spectral 
sequences (for any $r\geq 1$) between the cohomological spectral sequnence 
of $Y$ and the above homological spectral sequence of $(X,X-Y)$, which provides
exactly the ismorphism $\cap [X]:
Gr^W_{2n-t}H^{2n-s-t}(Y) \to Gr^W_{-t}H_{s+t}(X,X-Y)$.

e)\ $\Pi^{-1}_\s$ from III.3.1 induces the boundary operator
$H_*(X,X-Y,\Z)=H_*(U,\partial U,\Z)\to H_{*-1}(\partial U,\Z)$. 
\end{prop}
\begin{proof} \ The main result here is the degeneration of the spectral 
sequence at level 2, which follows via duality $d)$ from the {\em NCD} case.

Fix a triangulation $T$ and let $T'$ be its first barycentric subdivision. 
Then the  Poincar\'e duality map (cf. II.2.6) can be organized in the 
following morphism of double complexes.

Let $A_{st}^T(Y)=C_t^T(\tilde{Y}^{s+1})$, which form a double complex
 with $\partial $ and $i$, similarly as in III.1. Set the double complex
$B_{st}^T(Y)=Hom(A_{st}^T(Y),\Z)$ with dual morphisms $\delta$ and $i^*$.
Let $A_{s+1,t}^{T'}(X,X-Y)=C^{\p,T'}_{t+2s}(\tilde{Y}^{-s})$ with boundary 
morphisms $\partial$ and $\cap$ (similarly as $A_{**}(X,X-Y)$ defined above).
Then $pd:B_{-s,2n-t}^T(Y)\to A^{T'}_{st}(X,X-Y)$ satisfies
$\partial\circ pd=pd\circ \delta$ (cf. II.2.8) and $\cap\circ pd=pd\circ i^*$
(cf. II.2.9). 

Now,  $A_{**}(Y)$ is quasi-isomorphic to $A_{**}^T(Y)$ (i.e. their 
spectral sequences are the same for $r\geq 1$), and the later is dual to
$B_{**}^T(Y)$. Using II.2.9, $pd$ induces  an isomorphism at the level 
of the $E^1$ term, hence it is an quasi-isomorphism, and it induces
isomorphism at the level of any $E^r$ ($r\geq 2$). On the other hand,
$E^r(A_{**}^{T'}(X,X-Y))=E^r(A_{**}(X,X-Y))$ for $r\geq 1$.
Hence the result follows.
\end{proof}

\subsection{}{\bf The homological double complex of $X-Y$.} \\

\noindent Define 
$A_{s,t}(X-Y):=C_{t+2s}^\pitchfork (\tilde Y^{-s}),
$ with $s\leq 0$ and $t+2s\geq 0$. 
Then  $D=\partial + \cap$ is a 
differential of the total complex $\Tot_* (A_{**}(X-Y))$.

Similarly to the case of the pair $(X,Y)$, we can define $A_{**}^{\p}(X)$ 
by $A_{s*}^{\p}(X)= C^{\p}_*(X)$ if $s=0$ and $=0$ otherwise. Then 
$A_{s*}(X-Y)= A_{s*}^{\p}(X)\oplus A_{s+1,*}(X,X-Y)$. 
Actually, $A_{*-1,*}(X-Y)$ is the cone of the morphism
$\cap:A^{\p}_{**}(X) \to A_{**}(X,X-Y)$, where $\cap $ is the intersection for 
$s=0$ and zero otherwise. 

\vspace{2mm}

\noindent  Hence the Poincar\'e duality in III.3.2  extends to:

\subsubsection{ }\begin{prop}

a)\ $E_{s,t}^1=H_{t+2s}(\tilde Y^{-s})$, and 
$d_1$ is the transfer map $i_!$.

b)\ $E^{r}_{st}\Longrightarrow H_{s+t}(X-Y,\Z) $, and 
$E^{\infty}_{st}\otimes \Q=Gr^W_{-t}H_{s+t}(X-Y,\Q)$.

c)\ $ E^*_{**}\otimes \Q$ 
degenerates at level two, i.e. $d^r\otimes 1_{\Q}=0$ for $r\geq 2$. 

d)\ The Poincar\'e duality map $pd$ (cf. II.4.6)
 induces an isomorphism of spectral 
sequences (for any $r\geq 1$) between the cohomological spectral sequnence 
of $(X,Y)$ and the above homological spectral sequence of 
$X-Y$, which provides
exactly the ismorphism $\cap [X]:
Gr^W_{2n-t}H^{2n-s-t}(X,Y) \to Gr^W_{-t}H_{s+t}(X-Y)$.
\end{prop}

\subsubsection{}{\bf Remark.} Notice that by corollary II.3.6  and III.4.1
$$\Pi^{-1}_\s:(Tot_*(A_{**}(X-Y)),D) \to (C_*(Z),\partial)$$
is a quasi-isomorphism of complexes.

Notice also that $H_*(C_*(Z))=H_*(X-Y)$. Indeed, if $Z_{\alpha}$ is a small 
collar of $Z$, then $H_*(Z)=H_*(Z-Z_{\alpha})=H_*(Z-\partial Z)=H_*(X-Y)$.

\vspace{2mm}

\noindent The above morphism can be extended to the pair $(X,X-Y)$ as follows.
Consider now the complex $(Ker_{*+1},\partial)$ defined as the kernel of the 
morphism $\Pi_*:(C_*(\partial Z),\partial)\to (C_*(Y),\partial)$
induced by $\Pi$. In fact, this complex coincide with the kernel of the
morphism 
$ \Pi_*:(C_*(Z),\partial)\to (C_*(X),\partial)$
since $\Pi$ induces an isomorphism over $X-Y$. 

\subsubsection{}\begin{corollary} \ 

a)\ $\Pi^{-1}_\s:A_{s,t}(X-Y)\to C_{s+t}(Z)$
maps $A_{s+1,t}(X,X-Y)\subset A_{s,t}(X-Y)$ in $Ker_{s+t+1}$.

b)\ There is a natural commutative diagram 
whose vertical arrows are quasi-isomorphisms.
$$\begin{array}{ccccccccc}
0&\to& Tot_*(A_{*+1,*}(X,X-Y))& \to & Tot_*(A_{**}(X-Y))&\to
& Tot_*(A_{**}^{\p}(X))&\to & 0\\
&&\Big\downarrow\vcenter{\rlap{$\Pi^{-1}_\s$}}
&&\Big\downarrow\vcenter{\rlap{$\Pi^{-1}_\s$}}
&&\Big\downarrow\vcenter{\rlap{$j_1$}}&&\\
0&\to&Ker_{*+1}&\to&C_*(Z)&\to&C_*(X)&\to& 0\end{array}$$
The horizontal lines induce the 
long homology exact sequence of the pair $(X,X-Y)$.
\end{corollary}
\begin{proof} The image of any cycle  $\xi\in C^{\p}_k(\tilde{Y}^p)$, with $p\geq 1$, via 
the composite map
$$C^{\p}_k(\tilde{Y}^p)\stackrel{\Pi^{-1}_p}{\longrightarrow}
C_{k+p}(\partial Z)\stackrel{\Pi_*}{\longrightarrow}C_{k+p}(Y)$$
has $k$-dimensional support. Actually, this support is the image of $|\xi|$ by
the natural map $\tilde{Y}^p\to Y$. Hence it supports no cycle in
$C_{k+p}(Y)$. 
\end{proof}

\section{Review of the spectral sequence associated with a double complex}

\subsection{} \ 
Let $ A = \oplus_{s,t} A _{s,t} $ be a finite homological double
complex with operators
\begin{diagram}
& A_{s - 1, t} & \lTo^{\delta} & A_{s,t}    
  & \quad \delta^2 = 0,\, \d^2 = 0,\, \delta \, \d + \d \, \delta  = 0 \\
&              &             &\dTo_{\d} &            \\
&              &             & A_{s, t - 1} &            \\
\end{diagram}
Let $ \Tot_k A = \oplus_{ s + t = k} A _{s,t} $ and $ D = d + \delta $
be the associated total complex. Define the (weight) filtration by $
W_{\ell} A = \oplus_{s \leq \ell} A_{s, t} $. By general theory, there
is a spectral sequence $ (E^{r}_{s,t}, d_r)_{r \geq 1} $, converging
to the graded homology of $A$.
We put $ E^0_{s,t} = A_{s,t} $; then
$$
E^{1}_{s,t} = \Ker (\d: A_{s,t} \rightarrow A_{s, t - 1}) /Im
(\d : A_{s, t + 1} \rightarrow A_{s,t})
$$
and $ d^1 $ is induced at the homology level by $ \delta $. Let
$$
\eqalign{
Z^r_{s,t} & = \{ c \in W_s \Tot_{s + t} A: D c \in W_{s - r} A\}
\cr
Z^r_s     & = \{ c \in W_s A: D c \in W_{s - r} A \} \cr
Z^{\infty}_{s,t} & = \{ c \in W_s \Tot_{s + t} A: D c = 0 \}\cr
Z^{\infty}_s & = \{c \in W_s A : D c = 0\}.
}
$$

\noindent Then
$$
\eqalign{
E^r_{s,t} &= Z^r_{s,t} / (Z^{r - 1}_{s - 1} + D Z^{r - 1}_{s + r - 1})
\cr
E^\infty_{s,t} & = Z^\infty_{s,t} / (Z^\infty_{s - 1} + D A \cap W_s
A).\cr
}
$$

\noindent 
In particular, $ E^r_{s,t} $ is generated by the class of elements
$$
c_{st} + c_{{s - 1}, {t + 1}} + \cdots + c_{s - r + 1, t + r - 1} \in
Z^r_{st} \quad (c_{pq} \in A_{pq})
$$
which  satisfy
$$
\d c_{st} = 0; \quad \d c_{s - 1, t + 1} + \delta c_{st} = 0;\quad
\ldots;\quad \d c_{s - r+1, t + r-1} + \delta c_{s - r + 2, t + r - 2} = 0.
\leqno \hbox{$(*_r)$}
$$
The map $d^r$ can be defined as:
$$
d^r_{s,t} [c_{st} + c_{{s - 1}, {t + 1}} + \cdots + c_{s - r + 1, t + r - 1}]
 = [\delta  c_{s - r + 1, t + r - 1}] $$

For $ r = 1 $, this means that $ c_{st} \in A_{pq} $ with $ \d
c_{st} = 0 $. For $ r = 2 $, $Z^2_{st}$ is generated by
elements of type $ c_{st} + c_{s - 1, t + 1}$ such
that $\d c_{s,t} = 0, \d c_{s - 1, t + 1} + \delta c_{s,t} = 0$,
which generates $ ker d^1 $ on $ E^1 $.
For $ r = \infty $, $ E^\infty_{st} $ is generated by class of cycles
of type
$$
c_{st} + c_{s - 1, t + 1} + \cdots + c_{s - r, t + r} + \cdots
$$
satisfying
$$
\d c_{s t} = 0; \quad \d c_{s - r, t + r} +  \delta c_{s - r +
1, t + r - 1} = 0 \text{ for any } r \geq 1.
\leqno \hbox{$(*_{\infty})$}
$$

\subsubsection{}
Now, assume that $ d^2 = d^3 = \cdots = 0 $ which means $ E^2_{s t}
\equiv E^\infty_{s t} $. Then consider
$$
\Ker (d_1: E^1_{s t} \rightarrow E^1_{s - 1,t}) \rightarrow
E^\infty_{s t}.
$$
This shows that for any $ c^2 =  c_{s t} + c_{s - 1, t + 1} $
satisfying $(*_2)$, 
there is a class $ c^\infty = c'_{s t} + c'_{s - 1, t + 1}  + \cdots $
\ with $(*_{\infty})$,
 such that $ [c^2] = [c^\infty] $ in $ E^\infty_{s t} $. More
precisely, we show that we can choose the cycle $ c^\infty $ with $
c'_{s t} = c_{s t} $.

\subsubsection{}
\begin{prop}
\begin{enumerate}
\item[a)] Consider a class $ [c_{st}] \in \ker [d_1 : E^1_{s t}
\rightarrow E^1_{s - 1,t}] $ and fix a representative $ c_{s t} \in
A_{s t} $ with $ \d c_{s t} = 0 $ and $ \delta c_{s t} \in \d
A_{s - 1, t + 1} $ i.e., we can define a cycle $ c^2 = c_{s t} + c_{s
- 1, t + 1} \in Z^2_{s t} $. If $ d^2 = 0 $, then $ c^2 $ can be
replaced by $ \tilde{c}_2 = c_{s t} + \tilde{c}_{s - 1, t + 1} $ such
that $ [c^2] =
[\tilde{c}^2] $ in $ E^2_{s t} $ and $ \tilde{c}^2 $ can be completed
to a cycle $ \tilde{c}^3 = c_{s t} + \tilde{c}_{s - 1, t + 1}  +
\tilde{c}_{s - 2, t + 2} \in Z^3_{s t} $.

The fact that $ d^2 = 0 $ is equivalent to the fact that this can be
done for any cycle $ c^2 $.

\item[b)] More generally, consider an element 
$ c^r = \sum^ {r-1}_{i = 0} c_{s - i, t + i}
\in Z^{r}_{s, t} $ and assume that $ d^r = 0 $.  Then $ c^r $  can be replaced
by $ \tilde{c}^r = c_{s t} + \sum^{r-1}_{i = 1} \tilde{c}_{s - i, t + i}
\in Z^{r}_{s, t}  $ with $[c^r] = [ \tilde{c}^r] \in E^{r}_{s, t} $
such that $ \tilde{c}^r $ can be completed to $ \tilde{c}^{r + 1} = 
 \tilde{c}^r  +
\tilde{c}_{s - r, t + r} \in Z^{r + 1}_{s t}$.

Again, in the case of an arbitrary spectral sequence, if any $c_r$
can be completed to $\tilde{c}^{r+1}$ as above, then $d^r=0$. 
\item[c)]
In particular, if $d^2=d^3=\cdots=0$, then
\begin{enumerate}
\item[i)]
any representative $c_{st}$ (with $\partial c_{st}=0$ and $\delta
c_{st}\in\partial A_{s-1,t+1})$ of a class in $\Ker (d_1|E_{st}^1)$ can
be completed to a cycle
$$
c_{s,t}^\infty=c_{s,t}+\sum_{i\ge1} c_{s-i,t-i}\in Z_{st}^\infty;
$$
\item[ii)]
if $c_{st}^\infty$ and $c_{st}'{}^\infty$ are two liftings of $c_{st}$
 as in (i)
, then
$c_{st}^\infty-c_{st}'{}^\infty\in Z_{s-1}^\infty$;
\item[iii)]
if $c_{st}$ and $c_{st}'$ are representatives of the same class of
$\ker(d_1|E_{st}^1)$ (i.e.
$c_{st}-c_{st}'\in\d B_{s,t+1}$) and we complete $c_{st}$ and $c_{st}'$
 to
$c_{st}^\infty$ and $c_{st}'{}^\infty$, respectively, as in (i), then
$$
c_{st}^\infty-c_{st}'{}^\infty\in DA\cap  W_s+Z_{s-1}^\infty.
$$
\end{enumerate}
\end{enumerate}
\end{prop}

\section{Constructing rational cycles}

The aim of the present section is to use the additional degeneration 
property of the weight spectral sequences and the above review
in order to give a construction of  cycles which represent
$W_{-t}H_{s+t}(A,B)$, where $(A,B)=(Y,\emptyset),\ (X,Y),\ (X,X-Y),\ (X-Y,
\emptyset)$ (The case $(\partial U,\emptyset)$ is again separated 
in section VI.).
Moreover, we find topological properties and characterizations
of these cycles. We will keep the notations of the previous sections. 
All the homology groups are considered with rational coefficients.

\subsection{}{\bf Cycles in $(X,X-Y)$.}

\subsubsection{}\ First notice that $H_*(X,X-Y)=H_*(U,\partial U)$. In this 
subsection we will construct relative cycles in $(U,\partial U)$, 
i.e chains $\xi$ supported in $U$ with $|\partial \xi|\subset \partial U$.

For any $(s,t)$ with $s+t=k$ fix $c_{st}\in A_{st}(X,X-Y)$ with 
$\partial c_{st}=0$ and $\cap c_{st}\in im\partial$. Then it can be 
completed to a cycle $c^{\infty}_{s,t}=c_{st}+c_{s-1,t+1}+\cdots$
with $Dc^{\infty}_{s,t}=0$. Consider $\Pi^{-1}_\s c^{\infty}_{s,t}\in C_{k-1}
(\partial Z)$. Then by III.3.1 one has $\partial \Pi^{-1}_\s 
c^{\infty}_{s,t}=0$.
Fix a sub-analytic homeomorphism  $\phi_{\alpha}:\partial Z\times[0,\alpha]\to Z_{\alpha}$,
where $Z_{\alpha}$ is a collar of $\partial Z\subset Z$ (cf. I), and
write $U=\Pi\circ\phi_{\alpha}(\partial Z\times [0,\alpha])$. 
For any chain $\xi\in C_*(\partial Z)$, one can associate in a natural way
a new chain $\xi\times [0,\alpha]\in C_{*+1}(\partial Z\times [0,\alpha])$.
Then $\Pi_*\circ (\phi_{\alpha})_*(\Pi^{-1}_\s
c^{\infty}_{s,t}\times[0,\alpha]))
\in C_k(U)$ has boundary supported in $\partial U=\Pi\circ
\phi_{\alpha}(\partial
Z\times \{\alpha\})$. For simplicity, we will denote this relatice cycle
by $c_{s,t}^{rel}$. 

Notice that for $p\geq 1$, $|c_{s,t}^{rel}|\cap Y^p=|c_{s',k-s'}|$,
where $s':=\min\{s,-p+1\}$. More precisely, if $p\geq 1-s$, the intersection
$|c_{s,t}^{rel}|\cap Y^p$ is $|c_{-p+1,k+p-1}|$, which has dimension $k-p-1$.
For $p=-s$, the intersection 
$|c_{s,t}^{rel}|\cap Y^p$ is $|c_{s,t}|$, which has dimension $k-p-2$.

\subsubsection{}\begin{prop}

a)\ Fix  $k = s + t $. To any homology class $[c_{st}]  $  in
$H_{t+2(s-1)}(\tilde Y^{-s + 1})$ with  $ [c_{st}]
\cap \tilde Y^{-s + 2} = 0 $, the above construction gives a
relative cycle $c_{s,t}^{rel}$. Their classes  generate $W_{-t}H_k(X,X-Y)=
W_{-t}H_k(U,\partial U)$ (and they are well-defined modulo $W_{-t-1}H_k$). 

b)\ 
$W_{-t}H_k(X,X-Y)=\{[c]: \ \mbox{$c$ a relative cycle in $X$ 
with $|\partial c|\subset X-Y$ such that} $ \
$\dim|c|\cap Y^p=k-p-1\ \ \mbox{for}\ \
p\geq 1-k+t,
\ \mbox{and}\ <k-p-1\ \ \mbox{for}\ \ p=-k+t\}.$
In these intersection conditions we always assume that  $p\geq 1$, otherwise 
the assumption is considered empty. 
In particular, $W_{-k}H_k(X,X-Y)=H_k(X,X-Y)$. 

The above formula  is also equivalent to the following characterization:
$$W_{-t}H_k(X,X-Y)=\{[c]: \dim|c|\cap Y^{-k+t} <2k-t-1.\}$$
\end{prop}

\vspace{2mm}

This shows that the homeomorphism type of the pair $(X,Y)$
(or even of $(U,\partial U)$) determines
completely the weight filtration of $H_*(X,X-Y)$.

\subsubsection{}{\bf The support filtration.} \  Let $U^p\subset X$ be a 
small regular  neighbourhood of $Y^p$ in $X$ (cf. [R-S], ch. 3). Then for 
$U$ sufficiently small (with respect to $U^p$), and $p\geq 1$,
one can consider the group
$$S_pH_k(U,\partial U):= im\,\big(\, H_k(U^p,U^p\cap\partial U)\to 
H_k(U,\partial U)\,\big).$$
It is not difficult to see that the (decreasing) 
filtration $\{S_pH_k(U,\partial U)\}_p$ is independent on the choice of 
neighbourhoods $U^p$; 
it will be called the ``support filtration'' of $H_k(U,\partial U)$. 

In the  above construction, $|c_{s,t}^{rel}|\subset U^p$ (where $p=1-s$ and $s+t=k$),
therefore:
$$W_{-k-p+1}H_k(U,\partial U)\subset S_pH_k(U,\partial U).$$
In general, the inclusion is strict.

\subsection{Cycles in $X-Y$.}

\subsubsection{}\ Here we will use the identification
$H_*(X-Y)=H_*(Z)$. Similarly as above, 
for any $(s,t)$ with $s+t=k$ fix $c_{st}\in A_{st}(X-Y)$ with 
$\partial c_{st}=0$ and $\cap c_{st}\in im\partial$. Then it can be 
completed to a cycle $c^{\infty}_{s,t}=c_{st}+c_{s-1,t+1}+\cdots$
with $Dc^{\infty}_{s,t}=0$. Consider $\Pi^{-1}_\s c^{\infty}_{s,t}\in C_{k}
(Z)$ (cf. III.4.2).

\subsubsection{}\begin{prop}

a)\ The cycles $\Pi^{-1}_\s c^{\infty}_{s,t}$ generate  $W_{-t}H_k(Z)=
W_{-t}H_k(X-Y)$ (and they are well-defined modulo $W_{-t-1}$). 

b)\  $W_{-k}H_k(X-Y)=H_k(X-Y)$.  For $t<-k$, $W_{-t}H_k(X-Y)$
can be characterized by the exact sequence of the pair $(X,X-Y)$:
$$H_{k+1}(X)\to H_{k+1}(X,X-Y)\to H_k(X-Y)\to H_k(X).$$
Above $H_i(X)$ is pure of weight $-i$, and the weight filtration of
$H_{k+1}(X,X-Y)$  is characterized in the previous subsection.
\end{prop}

\vspace{2mm}

\noindent For a  more detailed and geometrical presentation of the cycles 
in $X-Y$, see theorem III.0.

\vspace{2mm}

The above result 
 shows that  the homeomorphism type of the pair $(X,Y)$ determines
completely the weight filtration of $H_*(X-Y)$. On the other hand, it is well 
known that one cannot recover the weight filtration from
the homeomorphism type (even the analytic type) of $X-Y$. For a counterexample,
see e.g. [S-S], (2.12).

\vspace{2mm}

\subsubsection{}{\bf The support filtration.} Similarly as above, let 
$U^p$ be a regular neighbourhood of $Y^p$ in $X$ for any $p\geq 1$, and 
for $p=0$ take $U^0=X$. Define for any $p\geq 0$:
$$S_pH_k(X-Y):=im\,\big(\, H_k(U^p\cap(X-Y))\to H_k(X-Y)\,\big).$$
This provides the decreasing ``support filtration'' $S_*H_k(X-Y)$. Then:
$$W_{-k-p}H_k(X-Y)\subset S_pH_k(X-Y),$$
and the inclusion, in general, is strict.

\subsection{}{\bf Cycles in $Y$.}

\subsubsection{}\ Consider $(A_{**}(Y),\partial,i)$ and fix a pair $(s,t)$
with $s+t=k$. 
If one wants to construct a closed cycle in $Y$, it is natural to
start with a chain  $c_{st}\in C_t(\tilde{Y}^{s+1})$ (with 
$\partial c_{st}=0$) and tries to extent it. The first obstruction is
$i(c_{st})\in im\partial$, i.e. the existence of $c_{s-1,t+1}$ with
$i(c_{st})+ \partial c_{s-1,t+1}=0$. 
Then the second obstruction is $i(c_{s-1,t+1})\in
im\partial$, and so one. The remarkable fact is that once the first 
obstruction  is satisfied then all the others are automatically satisfied.
This follows from the degeneration of the spectral sequence (III.1.2), 
and it is a 
consequence of the algebraicity of $Y$; in a simple topological context it 
is not true.

If $c_{st}\in A_{st}(Y)=E^0_{st}$ satisfies $\partial c_{st}=0$ and
$i(c_{st})\in im\partial$, then it can be completed to a cycle 
$$c^{\infty}_{st}=c_{st}+c_{s-1,t+1}+\cdots +c_{0,k}$$
with $Dc^{\infty}_{st}=0$ (see section IV). 
Let $n:\y\to Y$ be again the natural projection for any $p\geq 1$.
Then the correspondence 
$c_{st}\mapsto n_*(c_{0,k})\in C_*(Y)$ defines a closed cycle in $Y$. 
(We invite the reader to verify that $\partial n_*(c_{0,k})=0$.) 

\subsubsection{}{\bf Proposition.} {\em 
Fix $k=s+t$. To any homology class $[c_{st}]\in H_t(\tilde{Y}^{s+1})$ with 
$i[c_{st}]=0$ one can construct a cycle 
$$c^{\infty}_{st}=c_{st}+c_{s-1,t+1}+\cdots +c_{0,k}$$
with $Dc^{\infty}_{st}=0$. The homology classes $[n_*c_{0,k}]\in H_k(Y)$
generate $W_{-t}H_k(Y)$. Any other completion $c^{\infty}_{st}$ of $c_{st}$
provides the same homology class modulo $W_{-t-1}H_k(Y)$. }

\vspace{2mm}

Actually, the chains  $c_{s-i,t+i}$  can be recovered from the intersections
$|c_{0,k}|\cap Y^{s-i+1}$. Indeed,  for any $i$ with $s>i\geq 0$, and 
write $p=s-i+1$:
$$|n_*c_{0,k} |\cap Y^p=n(|c_{p-1,k-p+1}|),\ \mbox{whose dimension is}\ 
k-p+1.$$
On the other hand, if $p\geq s+2$, then 
in the homology class $[c_{st}]$  one can  take a representative $c_{st}$
which is in a general position with respect to the $Y^p$ (i.e. 
they are transversal), therefore:
$$|n_*c_{0,k} |\cap Y^p= |n_*c_{st} |\cap Y^p
\ \mbox{whose dimension is}\ t-2(p-s-1)<k-p+1.$$

\subsubsection{}\begin{corollary}
The weight filtration of $H_k(Y)$ can be characterized by:
$$W_{-t}H_k(Y)=\{[c]: \dim|c|\cap Y^p=k-p+1\ \mbox{for}\ \
p\leq k-t+1,$$
$$\hspace{3cm}\mbox{and}\ <k-p+1\ \mbox{for}\ \ p=k-t+2.\}$$
Actually, this is equivalent to the following characterization:
$$W_{-t}H_k(Y)=\{[c]: \dim|c|\cap Y^{k-t+2} < t-1.\}$$
\end{corollary}
This can be rewritten in the language of intersection homology as
follows. For any integer $s\geq 0$,
consider the perversity $\u{s}$ defined by $\u{s}(2i)=i$ for $0\leq i\leq s$,
and $\u{s}(2i)=s$ for $i\geq s$. (Since we have no stratum with odd
codimension, $\u{s}(2i+1)$ is unimportant.)
Notice that $\u{s}$ is not a perversity in the sense of [GM] (i.e. does not 
satisfy $\u{s}(2)=0$), it is a {\em generalized perversity} 
(see. e.g. [K] or [H-S]). 
\subsubsection{}{\bf Corollary.} \ \ 
$W_{-k+s}H_k(Y)=im\,(IH^{\u{s}}_k(Y)\to H_k(Y)).$

\vspace{2mm}

The above corollary in fact says that Deligne's weight  filtration and the 
Zeeman filtration (or ``support filtration'') coincide. This fact was 
conjectured by MacPherson, and verified by C. McCrory in [McC] (see also
[GNPP] and [H-S]). 
One of the main consequences of this fact is that the weight filtration
of $Y$ is completely topological. 

\vspace{2mm}

\noindent In the last part of this subsection, we present a different 
construction for the case $Y$. 

\subsubsection{}{\bf Milnor's construction.}\ 
Sometimes it is convenient to 
replace $Y$ by another space, 
called  the  ``geometric realization'' of the ``semi-simplicial object''
$\{\tilde{Y}^p\}_p$. 
The construction is done by Milnor [M] (see also [S], [D] and [A]).

This  new space is homotopic to $Y$. It can be constructed as follows. 
Consider the natural (normalization) map $n:\tilde{Y}^1\to Y$. For any points
$y\in Y$ with $y\in Y_{\alpha_1,\ldots,\alpha_j}$ and $y\not\in Y_{\beta}$
for $\beta\not\in \{\alpha_1,\cdots,\alpha_j\}$, the set $n^{-1}(y)$
consists of $j$ points $\{y_{k}\}_k$. Obviously, if for any $y$,
we identify in $\tilde{Y}^1$ the points $\{y_{k}\}_k$, then the quotient-space
is exactly $Y$. But we want a slightly different construction: 
for any $y$, we glue to $\tilde{Y}^1$ a $(j-1)$--simplex 
$\Delta_y=[y_1,\ldots,y_j]$, 
in a compatible way. The new space $SY$ will have the following properties. 
The map $sn:SN\to Y$, defined by $sn(\Delta_y)=y$ is continuous, it extends $n$
(i.e. $sn|\tilde{Y}^1=n$), and it is  a homotopy equivalence. 

Actually, $SY$ can be written as 
$$SY=\coprod_{p\geq 1} \tilde{Y}^p\times \Delta_{p-1}\ /\sim,$$
where $\Delta_{p-1}$ is 
the $(p-1)$--simplex, and $\sim$ is a suitable identification.
If for any $0\leq j\leq p-1$, $i_j:\tilde{Y}^p\to \tilde{Y}^{p-1}$
denote the natural maps induced by the inclusions
$Y_{\alpha_1,\ldots,\alpha_p}\hookrightarrow Y_{\alpha_1,\ldots,
\hat{\alpha}_{j+1},
\cdots,\alpha_p}$, and $\partial_j:\Delta_{p-2}\to \Delta_{p-1}$
are the face maps, then we identify:
$$(i_jy,x)\sim (y,\partial_jx), \ \mbox{for any}\ 
0\leq j\leq p-1, \ y\in\tilde{Y}^p,\ \mbox{and}\ \ x\in \Delta_{p-2}.$$

\noindent For more details,  see [M], [G] or [A] (with even 
some pictures on page 239).

\vspace{2mm}

One of the advantages of the space $SY$ is the existence of a morphism
$sn^{-1}:C_k(\tilde{Y}^p)\to C_{k+p-1}(SY)$.
This is defined as follows.  
For any $k$--simplex $\sigma$ we will denote 
by the same symbol $\sigma$ the chain with support
$\sigma$ and coefficient  one. We will  define 
$sn^{-1}(\sigma)$ by the chain 
$$[\sigma\times \Delta_{p-1}]\subset \coprod
\tilde{Y}^p\times \Delta_{p-1}\ /\sim.$$
(Here 
$[\sigma\times \Delta_{p-1}]$ can be covered by simplices, and we take all
of them with coefficient one; this is the wanted chain.)
This generates a linear map $sn^{-1}$. 

\subsubsection{}{\bf Lemma.} {\em The above map induces a morphism 
of complexes:
$$ns^{-1}:(Tot_*(A_{**}(Y)),D)\to (C_*(SY),\partial)$$
which is a quasi-isomorphism.}
\begin{proof} For the first part notice that 
$$
\partial(\sigma\times \Delta_{p-1})=
\partial \sigma \times \Delta_{p-1}+(-1)^k\sigma\times \sum_{i}(-1)^i
\partial_i\Delta_{p-1}=
\partial \sigma \times \Delta_{p-1}+i(\sigma)\times \Delta_{p-2},$$
i.e. $\partial (sn^{-1}(\sigma)=sn^{-1}(\partial \sigma)+sn^{-1}(i(\sigma))$
(cf. the identification $\sim$ above). The second part follows 
from  the above lemma III.1.2 (and its proof). 
\end{proof}

Then III.1.2 and IV.1.2 read as:
\subsubsection{}\begin{prop}

a)\ If $c_{st}\in A_{st}(Y)=E^0_{st}$ satisfies $\partial c_{st}=0$ and
$i(c_{st})\in im\partial$, then it can be completed to a cycle 
$$c^{\infty}_{st}=c_{st}+c_{s-1,t+1}+\cdots +c_{0,k}$$
with $Dc^{\infty}_{st}=0$. This means (cf. III.1.3) that $\partial 
ns^{-1}c^{\infty}_{st}=0$, generating a homology class in 
$H_k(SY)=H_k(Y)$. 

b)\ The homology classes $[ns^{-1}c^{\infty}_{st}]$ generate
$W_{-t}H_k(Y)$. Any other completion $c^{\infty}_{st}$ of $c_{st}$ provides
the same homology class modulo $W_{-t-1}H_k(Y)$. 
\end{prop}

The cycle $ns^{-1}(c^{\infty}_{st})$ constructed above
can be projected via $sn$ into $Y$. In this case, all the  chains 
$c_{s-i,t+i}$ (for $s>i\geq 0$) provide ``degenerate'' chains in $Y$
(with support of dimension $t+i$). This shows that the correspondence 
$c_{st}\mapsto sn^{-1}(c^{\infty}_{st})\in  C_*(SY)$ can be replaced by 
$c_{st}\mapsto n_*(c_{0,k})\in C_*(Y)$. 

\subsection{Cycles in $(X,Y)$.}

\subsubsection{}\ The construction is similar to the previous case. Take 
$c_{st}\in A_{st}(X,Y)$ 
with $\partial c_{st}=0$ and  $ i(c_{st})\in im \partial$.
Then complete to $c^{\infty}_{st}=c_{st}+\cdots +\tilde{c}_{0,k}$,
 where $k=s+t$.
Here $\tilde{c}_{0,k}\in C_k(X)$ with $\partial \tilde{c}_{0,k}
=-i(c_{1,k-1})$, hence
$|\partial \tilde{c}_{0,k}|\subset Y$ provided that $s>0$. 
\subsubsection{}\begin{prop}\ 

a)\ The relative cycles $\tilde{c}_{0,k}$ (associated with $c_{st}$ as above)
generate $W_{-t}H_k(X,Y)$, and are well-defined modulo $W_{-t-1}H_k(X,Y)$.

b)\ \ 
$W_{-t}H_k(X,Y)=\{[c]: \ \mbox{$c$ a relative cycle in $X$ with $|\partial c|
\subset Y$ such that} $ \
$\dim|c|\cap Y^p=k-p\ \ \mbox{for}\ \
p\leq k-t,
\ \mbox{and}\ <k-p\ \ \mbox{for}\ \ p=k-t+1\}.$
This is equivalent to the following characterization:
$$W_{-t}H_k(X,Y)=\{[c]: \dim|c|\cap Y^{k-t+1} < t-1\}.$$

c)\ If $\partial :H_k(X,Y)\to H_{k-1}(Y)$ is the natural map, then
$W_{-t}H_k(X,Y)=\partial^{-1}W_{-t}H_{k-1}(Y)$. In particular,
$W_*H_*(X,Y)$ is topological. Moreover, it can also be characterized by
$$W_{-t}H_k(X,Y)=\{[c]: \dim|\partial c|\cap Y^{k-t+1} < t-1\}.$$
\end{prop}

\subsubsection{} \ 
Let $i'$ be the composite of the natural maps 
$C_*(SY)\to C_*(Y)\to C_*(X)$
(induced by $sn$ and the inclusion). Then there is the following 
commutative diagram of complexes:
$$\begin{array}{ccccccccc}
0&\to&Tot_*(A_{**}(X))&\to &
Tot_*(A_{**}(X,Y))&\to &
Tot_{*}(A_{*-1,*}(Y))&\to&0\\
&&\Big\downarrow\vcenter{\rlap{$1$}} & &
\Big\downarrow\vcenter{\rlap{$1\oplus sn^{-1}$}} & &
\Big\downarrow\vcenter{\rlap{$sn^{-1}$ }}&& \\
0&\to&C_*(X)&\to& Cone_*(i')& \to & C_{*-1}(SY)&\to &0 \end{array}$$

It is clear that the second line provides the long homology exact sequence 
of the pair $(X,Y)$, and the vertical arrows in the diagram  are 
quasi-isomorphisms.

\subsection{The homology of $Y^p$ and $(X,X-Y^p)$.}

\subsubsection{} \ We end this section with the discussion of 
some  connections between
the weight filtration and the homology of the spaces $Y^p$ and
$(X,X-Y^p)$. The proofs are standard or similar to the proofs already 
presented, therefore we will omit them. 
If $A_{**}$ is a double complex, we denote by $\sigma_{s\leq i}A_{**}$ the
subcomplex of $A_{**}$ defined by $(\sigma_{s\leq i}A_{**})_{st}=A_{st}$
for $s\leq i$ and zero otherwise. $\sigma_{s\geq i+1}A_{**}$ is the quotient
double complex $A_{**}/\sigma_{s\leq i}A_{**}$. 

\subsubsection{}\ 
We start with the double complex of $Y$. For any $p\geq 1$:
$$H_k(Tot_*(\sigma_{s\geq p-1}\,A_{**}(Y)))=H_{k-p+1}(Y^p).$$
Moreover, $W_*(\sigma_{s\geq p-1}A_{**}(Y))$ induces the weight filtration 
on this space. The exact sequence:
$$0\to \sigma_{s\leq p-2}A_{**}(Y)\to A_{**}(Y)\to \sigma_{s\geq p-1}A_{**}(Y)
\to 0$$
provides the identity:
$$W_{p-2-k}H_k(Y)=\,\ker\,\big( \, H_k(Y)\stackrel{b}{\to}
 H_{k-p+1}(Y^p)\,\big),$$
where  $b$ is the (``Mayer-Vietoris'') boundary map (associated with the closed
covering $Y=\cup_iY_i$).

\vspace{2mm}

Similar identities are valid for $A_{**}(X,Y)$, in particular, for $p\geq 1$:
$$W_{p-1-k}H_k(X,Y)=\,\ker\,\big( \, H_k(X,Y)\stackrel{b'}{\to}
 H_{k-p}(Y^p)\,\big),$$
where  $b'$ is the composite of the boundary operator and $b$.

\vspace{2mm}

In the case of $(X,X-Y)$, for $p\geq 1$, one has:
$$H_k(Tot_*(\sigma_{s\leq -p+1}\,A_{**}(X,X-Y)))=H_{k+p-1}(X,X-Y^p),$$
and  $W_*(\sigma_{s\leq -p+1}A_{**}(X,X-Y))$ induces the weight filtration 
on this space. The exact sequence:
$$0\to \sigma_{s\leq -p+1}A_{**}(X,X-Y)\to A_{**}(X,X-Y)\to 
\sigma_{s\geq -p+2}A_{**}(X,X-Y) \to 0$$
provides:
$$W_{-p+1-k}H_k(X,X-Y)=\,im\,\big( \, H_{k+p-1}(X,X-Y^p)\stackrel{b}{\to}
 H_{k}(X,X-Y)\,\big),$$
where  $b$ again is a  (``Mayer-Vietoris'') boundary map
(associated with the open covering $(X,X-Y_i)_i$).

\vspace{2mm}

\noindent Similar property holds  for $A_{**}(X-Y)$, in particular, for
$p\geq 1$:
$$W_{-p-k}H_k(X-Y)=\,im\,\big( \, H_{k+p}(X,X-Y^p)\stackrel{b'}{\to}
 H_{k}(X-Y)\,\big),$$
where  $b'$ is the composite of $b$ and the boundary map. 

\vspace{2mm}

The above properties of the pairs $(Y,\emptyset)$ and $(X,X-Y)$, respectively of
$(X,Y)$ and $(X-Y,\emptyset)$ correspond by the Poincar\'e Duality. 
This follows also from the commutativity of some diagrams like the following:

$$\begin{array}{ccc}
H^{2n-k}(Y)&\stackrel{PD}{\longrightarrow}&
H_{k}(X,X-Y)\\
\Big\uparrow\vcenter{\rlap{$b^*$}}&&
\Big\uparrow\vcenter{\rlap{$b$}}\\
H^{2n-k-p+1}(Y^p)&\stackrel{PD}{\longrightarrow}&
H_{k+p-1}(X,X-Y^p)\end{array}$$
Above, $PD$ denotes the Poincar\'e Duality \ $\cap[X]$ (induced by $pd$). 

\section{The homology   of $\partial U$.}

\subsection{}{\bf The homological double complex of $\partial U$.} 

\subsubsection{}\ 
The space $\partial U$ appears in a natural way in two homological exact sequences. One of 
them is the pair
$(U,\partial U)$. The homological exact sequence of this pair, and the
isomorphism $H_*(U)=H_*(Y)$ suggests that a possible double complex for 
$\partial U$  should satisfy 
$$0\to A_{*+1,*}(U,\partial U)\to A_{**}(\partial U)\to A_{**}(Y)\to 0.$$
Here, the double complexes of $(U,\partial U)$, by definition is the same 
as the double complex of $(X,X-Y)$. Now, using this and the double complex 
of  $Y$, one can try to define the double
complex of $\partial U$ by
$$A_{s,t}(\partial U)=\left\{\begin{array}{ll}
C^{\p}_{t+2s}(\tilde{Y}^{-s})&\  \mbox{for}\ s\leq -1\\
C_t(\tilde{Y}^{s+1})&\ \mbox{for}\ s\geq 0.\end{array}\right.$$
But now we are confronted with the definition of the arrows of the double
complex. Actually, we have all the vertical arrows, and all the horizontal
 arrows corresponding to $s\leq -1$ and $s\geq 0$. But we need to define 
a map $C_*(\tilde{Y}^1)\to C_{*-2}^{\p}(\tilde{Y^1})$ with some nice
properties.
First of all, this map should be compatible with the other arrows, in the
sense that the whole complex should form a double complex. On the other hand, 
we expect (since we know the cohomological $E_1$ term of the 
corresponding spectral sequence)
that this map should induce at the homology level  the ``intersection matrix'';
more precisely, the map $\oplus_{\alpha}c_{\alpha}\mapsto \oplus_{\alpha}
d_{\alpha}$, where for $c_{\alpha}\in H_k(Y_{\alpha})$ one gets  $d_{\alpha}=
\sum_{\beta}\,  c_{\beta}\cap [Y_{\alpha}]$.
Here, if $\alpha\not=
\beta$, then $c_{\beta}\cap[Y_{\alpha}]$ is provided by the transfer map 
$i_!:H_k(Y_{\beta})\to H_{k-2}(Y_{\alpha,\beta})$ composed with
$i_*:H_{k-2}(Y_{\alpha,\beta})\to H_{k-2}(Y_{\alpha})$ induced by the 
inclusion. If $\alpha=\beta$, then $c_{\alpha}\cap [Y_{\alpha}]$ is the
cap product by $[Y_{\alpha}]$. 

Therefore, the wanted 
map $C_*(\tilde{Y}^1)\to C_{*-2}^{\p}(\tilde{Y^1})$ should be some kind of intersection, but 
we realize immediatelly that we face  serious obstructions:
we have to intersect cycles which are not ``transversal'', and the image should be special 
``transversal'' cycle. Even if we try to modify our complexes,
similar type of obstruction will survive. The explanation for this is 
the following. The cap product 
$c_{\alpha}\cap [Y_{\alpha}]$ cannot be determined only from the 
spaces $\{\tilde{Y}^p\}_{p\geq 1}$, we need the Poincar\'e dual of 
the spaces $Y_{\alpha}$ in $X$, hence we need also to consider  the space $X$
(or at least $U$) in our double complex.

So, we have to think about $\partial U$ as the boundary of $Z$, and we have to consider the 
pair of spaces $(Z,\partial Z)$. Notice that $H_*(Z,\partial Z)=
X_*(X,Y)$ and $H_*(Z)=H_*(X-Y)$, so it is natural to ask for for a double complex with:
$$0\to A_{*+1,*}(X,Y)\to A_{**}(\partial U)\to A_{**}(X-Y)\to 0.$$
The construction is done in the next subsection.

\subsubsection{}{\bf The construction of $A_{**}(\partial U)$.} \ 
Consider the following  diagram of complexes:
$$\begin{array}{ccccccccccc}
\cdots&
C_{*-4}^{\p}(\tilde{Y}^2)&\stackrel{\cap}{\leftarrow}&
C_{*-2}^{\p}(\tilde{Y}^1)&\stackrel{\cap}{\leftarrow}&
C_{*}^{\p}(\tilde{Y}^0)& &&&& \\
&&&&& \Big\downarrow\vcenter{\rlap{$j_1$}}&&&&& \\
&&&&& 
C_*(\tilde{Y}^0)& \stackrel{i}{\leftarrow}&
C_*(\tilde{Y}^1)& \stackrel{i}{\leftarrow}&
C_*(\tilde{Y}^2)& \cdots \end{array}$$
The first line corresponds to the double complex $A_{**}(X-Y)$ (where the
 column $s=0$ is $C_*^{\p}(\tilde{Y}^0)=C_*^{\p}(X)$), and the second line is the 
double complex $A_{**}(X,Y)$. Notice that $j_1$ can be considered as a 
morphism of double complexes, hence the usual cone construction 
provides  the double complex $A_{**}(\partial U)=\{A_{s*}(\partial U)\}_{s}$:
$$\begin{array}{ccccccccc}
\cdots&
C_{*-4}^{\p}(\tilde{Y}^2)&\stackrel{\cap}{\leftarrow}&
C_{*-2}^{\p}(\tilde{Y}^1)&\stackrel{\cap}{\leftarrow}&
C_{*}^{\p}(\tilde{Y}^0) &&& \\
&&& \oplus &  \stackrel{j_1}{\swarrow} &\oplus &&& \\
&&& 
C_*(\tilde{Y}^0)& \stackrel{i}{\leftarrow}&
C_*(\tilde{Y}^1)& \stackrel{i}{\leftarrow}&
C_*(\tilde{Y}^2)& \cdots \\
&&&&&&&&\\
& s=-2 && s=-1&&s=0&& s=1 & \end{array}$$

\subsubsection{}\begin{prop}

a)\ The $E^1$ term of the spectral sequence associated with 
$(A_{**}(\partial U),W)$ is:
$$\begin{array}{ccccccccc}
\cdots&
H_{*-4}(\tilde{Y}^2)&\stackrel{\cap}{\leftarrow}&
H_{*-2}(\tilde{Y}^1)&\stackrel{\cap}{\leftarrow}&
H_{*}(\tilde{Y}^0) &&& \\
&&& \oplus &  \stackrel{id}{\swarrow} &\oplus &&& \\
&&& 
H_*(\tilde{Y}^0)& \stackrel{i_*}{\leftarrow}&
H_*(\tilde{Y}^1)& \stackrel{i_*}{\leftarrow}&
H_*(\tilde{Y}^2)& \cdots \\
&&&&&&&&\\
& s=-2 && s=-1&&s=0&& s=1 & \end{array}$$

\noindent where  $\cap$ is given in II.2.5, e.g. 
$d\in H_k(X)$:
$$\cap(d)=\oplus_{\beta}\, (-1)^kd\cap [Y_{\beta}]\in 
\oplus_{\beta}H_{k-2}(Y_{\beta})=H_{k-2}(\tilde{Y}^1);$$
and $i_*$ is induced by $i$ (cf. III.1),  e.g. 
for $\oplus_{\alpha}
c_{\alpha}\in  \oplus_{\alpha} H_k(Y_{\alpha})=H_k(\tilde{Y}^1)$:
$$i_*(\oplus_{\alpha} c_{\alpha})=\sum _{\alpha}(-1)^k\, i_{\alpha}(c_{\alpha})\, \in 
H_k(X).$$

\noindent 
The above $E^1$ term is quasi-isomorphic to the complex:
$$\cdots H_{*-4}^{\p}(\tilde{Y}^2)\stackrel{\cap}{\leftarrow}
H_{*-2}^{\p}(\tilde{Y}^1)\stackrel{I}{\leftarrow}
H_*(\tilde{Y}^1) \stackrel{i_*}{\leftarrow}
H_*(\tilde{Y}^2) \cdots,$$
where $I:=\cap\circ i_*$, i.e.
$$I(\oplus_{\alpha}c_{\alpha})=\oplus_{\beta}\, (\sum_{\alpha}\, i_{\alpha}(
c_{\alpha})\cap [Y_{\beta}]\, ).$$

b)\ $E^r_{st}\Longrightarrow H_{s+t}(\partial U,\Z)$ and 
$E^{\infty}_{st}\otimes \Q=Gr^W_{-t}H_{s+t}(\partial U,\Q)$.

c)\ $ E^*_{**}\otimes \Q$ 
degenerates at level two, i.e. $d^r\otimes 1_{\Q}=0$ for $r\geq 2$. 
\end{prop}
\begin{proof} For $a)$ and $b)$ 
use  II.2.1 and II.2.5 and the corresponding definitions.
 $c)$  follows from the results of the previous subsections
and from the construction of the mixed cone.
\end{proof}

Now, we would like  to construct a morphism  of complexes $Tot_*(A_{**}
(\partial U))\to C_*(\partial U)$. For this the above complex is not
convenient, because of the presence of the global terms
 $C_*^{\p}(\tilde{Y}^0)$ and $C_*(\tilde{Y}^0)$. In the next construction,
we replace these complexes by the complexes of chains  supported by 
the close neighbourhood $U$. More precisely, we define $A'_{**}(\partial U)$:
$$\begin{array}{ccccccccc}
\cdots&
C_{*-4}^{\p}(\tilde{Y}^2)&\stackrel{\cap}{\leftarrow}&
C_{*-2}^{\p}(\tilde{Y}^1)&\stackrel{\cap}{\leftarrow}&
C_{*}^{\p}(U) &&& \\
&&& \oplus &  \stackrel{j_1}{\swarrow} &\oplus &&& \\
&&& 
C_*(U)& \stackrel{i}{\leftarrow}&
C_*(\tilde{Y}^1)& \stackrel{i}{\leftarrow}&
C_*(\tilde{Y}^2)& \cdots  \end{array}$$
Then  $E^1(A'_{**})$ is:
$$\begin{array}{ccccccccc}
\cdots&
H_{*-4}(\tilde{Y}^2)&\stackrel{\cap}{\leftarrow}&
H_{*-2}(\tilde{Y}^1)&\stackrel{\cap}{\leftarrow}&
H_{*}(U) &&& \\
&&& \oplus &  \stackrel{id}{\swarrow} &\oplus &&& \\
&&& 
H_*(U)& \stackrel{i_*}{\leftarrow}&
H_*(\tilde{Y}^1)& \stackrel{i_*}{\leftarrow}&
H_*(\tilde{Y}^2)& \cdots  \end{array}$$
which is quasi-isomorphic to $E^1(A_{**})$, hence $E^r(A_{**})=E^r(A'_{**})$
for any $r\geq 2$. 

\subsubsection{} \ 
Now, the advantage of this second double complex lies in the following 
construction.  Let $A_{**}(Z_{\alpha})$ be the double complex:

$$\begin{array}{ccccccccc}
\cdots&
C_{*-4}^{\p}(\tilde{Y}^2)&\stackrel{\cap}{\leftarrow}&
C_{*-2}^{\p}(\tilde{Y}^1)&\stackrel{\cap}{\leftarrow}&
C_{*}^{\p}(U) &\leftarrow &0&\cdots \\
&&&&&&&&\\
& s=-2 && s=-1&&s=0&&s=1  & \end{array}$$
Then, there is a natural morphisms of complexes $\Pi^{-1}_\s:
Tot_*(A_{**}(Z_{\alpha}))\to C_*(Z_{\alpha})$ (cf. II.3.8). This can be 
inserted in the following commutative diagram:
$$\begin{array}{ccccccccc}
0&\to& Tot_*(A_{*+1,*}(X,X-Y))& \to & Tot_*(A_{**}(Z_{\alpha}))&\to
& C_*^{\p}(U)&\to & 0\\
&&\Big\downarrow\vcenter{\rlap{$\Pi^{-1}_\s$}}
&&\Big\downarrow\vcenter{\rlap{$\Pi^{-1}_\s$}}
&&\Big\downarrow\vcenter{\rlap{$j_1$}}&&\\
0&\to&Ker_{*+1}&\to&C_*(Z_{\alpha})&\to&C_*(U)&\to& 0\end{array}$$

\vspace{2mm}

Since, using a  sub-analytic homeomorphism
 $Z_{\alpha}\approx \partial U\times [0,\alpha]$,
the second line in the above diagram induces the 
long homology exact sequence of the pair $(U,\partial U)$.

\subsubsection{}{\bf Corollary.}
{\em The vertical arrows in the above diagram are quasi-isomorphisms.}

\subsection{}{\bf Cycles in $\partial U$.}

\subsubsection{}\  
The complex $A_{**}(\partial U)$ was important from theoretical 
point of view: its spectral sequence is degenerating at rank two. 
On the other hand, the complex $A_{**}'(\partial U)$ has the same $E^r$ term
for $r\geq 2$, hence the degeneration property is still valid, but has the
advantage that it is local: $X$ is replaced by the neighbourhood $U$ of $Y$.
In the above construction we will use this second double complex $A_{**}'$.

Notice that $\sigma_{s\leq -1}A_{**}(X-Y)=
A_{*+1,*}(X,X-Y)$ is a subcomplex of $A_{**}'(\partial  U)$. 
The corresponding pair provides the exact sequence:

\subsubsection{}\hspace{1cm}
$\to H_{k+1}(X,X-Y)\stackrel{\partial }{\longrightarrow}
H_k(\partial U) \stackrel{\alpha }{\longrightarrow}H_k(Y)
\stackrel{\delta}{\longrightarrow} H_k(X,X-Y)\to \ .$

\vspace{2mm}

Since the weights of $H_{k+1}(X,X-Y)$ are $\leq -k-1$, and the weights of
$H_k(Y)$ are $\geq -k$, the weight filtration of $H_k(\partial U)$
is completely determined from the above exact sequence and from the 
weight filtrations of $(X,X-Y)$ and $Y$. For $l\leq -k-1$: $W_lH_k
(\partial U)=\partial (W_lH_{k+1}(X,X-Y)$, and for $l\geq -k$:
$W_lH_k(\partial U)=\alpha^{-1}W_lH_k(Y)$.

\vspace{2mm}

Notice again that the weight filtration of $H_k(\partial U)$ is completely
determined from the topology of the pair $(X,Y)$. On the other hand, it is
impossible to determine it from the diffeomorphism type of $\partial U$.
For counterexamples, see  again [S-S], e.g. (3.1). (For this one has to
use the fact that the link of on isolated singularity $(S,s)$
can be identified with 
the boundary $\partial U$ , where $U\to S$ is a  resolution  of 
the singular point  and $S$ is a Stein representative of $(S,s)$.)

\vspace{2mm}

Corresponding to the above discussion of the weights, 
in the construction of the cycles in $\partial U$
we will distinguish two different cases as well. 

\vspace{2mm}

Fix a pair $(s,t)$ with $s+t=k$ and $s\leq -1$. Consider $c_{st}\in
A_{s+1,t}(X,X-Y)\subset A_{s,t}(\partial U)$ with $\partial c_{st}=0$ and
$\cap c_{st}\in im \partial $. Complete to $c^{\infty}_{s,t}$ with
$Dc^{\infty}_{s,t}=0$. Then (cf. III.3.1) $\Pi^{-1}_\s c^{\infty}_{s,t}\in 
C_k(\partial Z)$ is closed. Recall that we have a natural identification
of $\partial Z$ and $\partial U$.

\subsubsection{}\begin{prop}
For $s\leq -1$, the cycles $\Pi^{-1}_\s c^{\infty}_{s,t}$ generate 
$W_{-t}H_k(\partial U)$. 
\end{prop}

Moreover, $[\Pi^{-1}_\s c^{\infty}_{s,t}]=\partial [c_{s,t}^{rel}]$,
where $c_{s,t}^{rel}$ is constructed in V.1. In particular, these homology 
classes inherit all the properties of $W_*H_{k+1}(X,X-Y)$
(including their relationship with the ``support'' filtration).
The details are left to the reader. 

\vspace{2mm}

Now, assume that $s\geq 0$. In this case the construction of the cycles
is more involved: we have to lift some cycles from $Y$ to $\partial U$.

Assume  that $c_{st}'\in A_{st}'(\partial U)$ ($s\geq 0$) satisfies 
$\partial c_{st}'=0$ and $d_1[c_{st}']=0$. If $s=0$ this means that 
$c_{st}'=c_{0k}'=c_{0,k}^{\perp}+c_{0,k}$, where $c_{0,k}^{\perp}\in
C_k^{\p}(U)$ and $c_{0,k}\in C_k(\tilde{Y}^1)$ such that $\partial 
c_{0,k}^{\perp}=\partial c_{0,k}=0$ and $i_*(c_{0,k})+c_{0,k}^{\perp}+
\partial \gamma=0$ for some $\gamma\in C_{k+1}(U)$. 

If $s>0$, then $c_{st}'\in A_{st}'(\partial U)$ has ``only one component'':
$c_{st}'=c_{st}\in C_t(\tilde{Y}^{s+1})$ with $\partial c_{st}=0$ and
$i_*c_{st}\in im\partial$. 

In both cases $c_{st}'$  can be completed to a cycle 
$$c^{\infty}_{s,t}=c_{s,t}+c_{s-1,t+1}+\cdots+c_{1,k-1}+(c_{0,k}+
c_{0,k}^{\perp})+
(\gamma_{-1,k+1}+c_{-1,k+1})+c_{-2,k+2}+\cdots$$
where $c_{s,t}\in C_t(\tilde{Y}^{s+1})$ for $s\geq 0$, $c_{s,t}\in 
C_{t+2s}^{\p}(\tilde{Y}^{-s})$ for $s\leq -1$, $c_{0,k}^{\perp}\in C_k^{\p}(U)$,
and $\gamma_{-1,k+1}\in C_{k+1}(U)$. These chains satisfies the following
relations:

\vspace{2mm}

$\partial c_{s,t}=0$;

$i_*c_{s-l,t+l}+\partial c_{s-l-1,t+l+1}=0$ for $0\leq l\leq s-1$;

$\partial c_{0,k}^{\perp}=0$;

$i_*c_{0,k}+jc_{0,k}^{\perp}+\partial \gamma_{-1,k+1}=0$;

$\cap c_{0,k}^{\perp}+\partial c_{-1,k+1}=0$;

$\cap c_{l,k-l}+\partial c_{l-1,k-l+1}=0$ for $l\leq -1$.

\vspace{2mm}

\noindent Now, we will separate the chain (for the notations, see  VI.1.4):
$$\bar{c}_{st}:=c_{0,k}^{\perp}+c_{-1,k+1}+c_{-2,k+2}+\cdots\ \in A_{**}(
Z_{\alpha}).$$
Then by the above relations, $D(\bar{c}_{st})=0$, where here $D$ is the 
differential in $A_{**}(Z_{\alpha})$. Therefore, $\Pi^{-1}_\s\bar{c}_{st}
\in C_k(Z_{\alpha})$ is closed. Obviously, the projection 
$pr:Z_{\alpha}\approx \partial Z\times [0,\alpha]\to \partial Z$
provides a closed cycle $pr_*\Pi^{-1}_\s \bar{c}_{st}\in C_k(\partial Z)$. 

\subsubsection{}\begin{prop}
For any $s\geq 0$ and $t+s=k$, the cycles
$pr_*\Pi^{-1}_\s \bar{c}_{st}$ generate $W_{-t}H_k(\partial U)$. \end{prop}

Notice that in the cycles $\bar{c}_{st}$ we do not see the chains 
$c_{l,k-l}$ for $l\geq 0$, in particular neither $c_{s,t}$, the chain which
generates $\bar{c}_{st}$. The chain $c_{st}$ is completed to a sequence
$c_{st}+c_{s-1,t+1}+\cdots + c_{0,k}$ with $c_{0,k}\in C_k(\tilde{Y}^1)$. The 
chain $i_*c_{0,k}$ actually is closed (in $U$) and supported by $Y$.
Now, this  is replaced by the transversal chain $c_{0,k}^{\perp}
\in C_k^{\p}(U)$ such that $i_*c_{0,k}+c_{0,k}^{\perp}+\partial \gamma=0$
for some $\gamma$, and finally $c_{0,k}^{\perp}$ is completed to
$\bar{c}_{st}$. It is really remarkable that the above algebraic construction
plays the role of a very geometric operation: it replaces a closed
chain supported in 
$Y$ by another closed chain supported in $U$, homologous with the original one
in $U$, and dimensionally
transversal to the stratification given by $Y$ (i.e. it solves the 
problem, obstruction mentioned in VI.1.1.)

\subsubsection{}{\bf Example.} Assume that $n=2$, and $Y$ is a connected
set of curves $\{Y_i\}_i$ in $X$. Set $g=\sum_i\mbox{genus}(Y_i)$.
If $\Gamma$ is the dual graph of the curves then let $c_{\Gamma}=rank 
H_1(|\Gamma|)$ be the number of independent cycles in $|\Gamma|$. 
Let $I$ be the intersection matrix of the irreducible curves $Y_i$. Then it is well-known that 
$rank H_1(\partial U)=rank\ker I+2g+c_{\Gamma}$. 
These three contributions correspond exactly to the weight of $H_1(\partial
U)$. 

Indeed, for $s=-1$, take $c_{-1,2}\in C_0^{\p}(\tilde{Y}^1)$ as above. A 
possible chain $c_{-1,2}$ is an arbitrary point $P$ in $Y^1-Y^2$
 with coefficient one. 
Then $c_{-1,2}^{\infty}=c_{-1,2}$ and $\Pi^{-1}(c_{-1,2})$ is a circle $S^1$,
the loop around $Y$ in a transversal slice at $P$. These loops $\gamma_P$
generate $W_{-2}H_1$ (isomorphic  to $coker\, I$). 
 If $s=0$, consider a closed 1-cycle $c_{0,1}$ in one of the components of $Y$.
This can be changed by a homologous cycle $c_{0,1}^{\perp}$ in $U$, which has
 no intersection with $Y$. Hence it can be contracted to $\partial U$. These 
cycles generate $W_{-1}H_1$ (so that $\dim Gr_{-1}^WH_1=2g$). 
Notice that the lifting
$c_{0,1}\mapsto c_{0,1}^{\perp}$ is defined modulo the cycles of type
$\gamma_P$. Finally, consider $c_{1,0}\in C_0(Y^2)$ such that $d_1[c_{1,0}]=0$.
Here $d_1:H_0(Y^2)\to H_0(\tilde{Y}^1)$ and $\ker d_1\approx H_1(|\Gamma|)$.
Take $c_{0,1}\in C_1(\tilde{Y}^1)$ such that $ic_{1,0}+\partial c_{0,1}=0$.
Then $i(c_{0,1})$ is a cycle in $U$ supported by $Y$. We replace it by 
$c_{0,1}^{\perp}$  which has no intersection with $Y$. They provide the
remaining cycles in $W_0H_1$ so that  $\dim Gr^W_0H_1=c_{\Gamma}$. 

Now, we discuss the case $H_2(\partial U)$ as well. Take  a point $P\in 
Y^2$ (with coefficient one) corresponding to $c_{-2,4}$. 
Then $\Pi^{-1}c_{-2,4}$ is a torus in $\partial 
U$. They generate $W_{-4}H_2$. Notice that $W_{-4}H_2\approx
coker \cap:H_2(\tilde{Y}^1)\to
H_0(Y^2)$ has dimension $c_{\Gamma}$. 
Next, take a generic closed 1-cycle in $\tilde{Y}^1$ whose image $c$ in $Y$
has no intersection with $Y^2$. The 2-cycle $\Pi^{-1}(c)$ in $\partial U$
is an $S^1$ bundle over $c$; they generate $W_{-3}H_2$. Finally, consider
$c_{02}=\sum_i m_iY_i\in C_2(\tilde{Y}^1)$ such that $\cap [c_{02}]=0$.
This means that $[c_{02}]\cdot [Y_j]=0$ for all $j$. The dimension of the space
generated by these classes $[c_{02}]$ is $coker\,I$.
Now, each $c_{02}$ is replaced by a transversal 2-cycle $c_{02}^{\perp}$.
Transversality implies that $c_{02}^{\perp}\cap Y_i$ is a 0-cycle in $Y_i$,
which by the above assumption is zero-homologous. In particular, 
$c_{02}^{\perp}\cap Y_i=\partial c_{-1,3}$. Now, take a very small tubular
neighbourhood $U'$ of $Y^1-Y^2$ with projection $pr:U'\to Y^1-Y^2$. 
(Here $U-U'$ stays for $Z_{\alpha}$.) Then 
the boundaries of the chains
$c_{02}^{\perp}-U'$ and the $S^1$-bundle $(pr|\partial U')^{-1}(c_{-1,3})$ can be identified 
(modulo sign) so they can be glued. They
 generate the remaining classes in $W_{-2}H_2$. 

\subsubsection{}{\bf Remark. -- Purity results.} \ The 
above example shows that 
in general all the possible weights (permitted by the spectral sequence)
can appear. For example, if $\partial U$ is the boundary (or link) of a 
1-parameter family of projective curves over a small disc, then the
intersection matrix $I$ has 1-dimensional kernel, hence $H_1(\partial U)$
can have weight $-2,-1$ and 0. On the other hand, if $\partial U$ is the 
boundary  of the resolution of a normal surface singularity, then $I$
 is non-degenerate, hence the non-trivial weight of $H_1(\partial U)$ 
are $-1$ and 0. 

More generally, if $S$ is a projective algebraic variety with unique 
singular point $s\in S$, and $\phi:X\to S$ is a resolution of this isolated 
singularity with normal crossing exceptional divisor 
$Y=\phi^{-1}(s)$, then the following additional restrictions
hold. They can be easily deduced from the  corresponding cohomological 
statements about the weight filtrations of the links of isolated singularities.
For a complete proof and details see  e.g. [N-A].

\vspace{2mm}

1) \ If $k\leq n-1$ then $Gr^W_lH_k(\partial U)=0$ for $l$ not in $[-k,0]$;

2)\ \ If $k\geq n$ then $Gr^W_lH_k(\partial U)=0$ for $l$ not in $[-2n,-k-1]$;

3)\ If $k\geq n$ then $H_k(Y)$ is pure of weight $-k$ and
$H_k(U)\to H_k(U,\partial U)$ is injective;

4)\ If $k\leq n$ then $H_k(U,\partial U)$ is pure of weight $-k$ and
$H_k(U)\to H_k(U,\partial U)$ is surjective;

5)\ Consider the exact sequence (cf. the $E^1$ term in VI.1.3):
$$\cdots H_{k-4}^{\p}(\tilde{Y}^2)\stackrel{\cap}{\leftarrow}
H_{k-2}^{\p}(\tilde{Y}^1)\stackrel{I_k}{\leftarrow}
H_k(\tilde{Y}^1) \stackrel{i_*}{\leftarrow}
H_k(\tilde{Y}^2) \cdots.$$
If $k\geq n$ then $\ker I_k= im\, i_*$, if $k\leq n$ then $\ker \cap = im\, 
I_k$. 

Actually, between the above exact sequence and the the exact sequence VI.2.2
there is the following connection: $\delta$ is non-trivial only for weight 
$-k$ and $Gr^W_{-k}\delta $ can be identified 
with $\hat{I}_k:H_k(\tilde{Y}^1)/im\, i_*\to \ker\cap$.

\end{document}